\documentclass[english,11pt,article]{amsart}
\usepackage[english]{babel}
\usepackage{amsmath,amssymb,enumerate,amsthm}
\usepackage[latin1]{inputenc}
\usepackage[T1]{fontenc}
\usepackage{cite}
\usepackage{latexsym}
\usepackage{a4}

\usepackage{psfig,labelfig}
\usepackage{graphicx}
\usepackage[dvips]{epsfig}

\newcommand{\co}{\colon\thinspace}   
\newcommand{\fnote}[1]{\footnote{\small sharp1}}
\newcommand{\vol}{\mbox{vol}}
\newcommand{\emvol}{\mbox{\em vol}}

\newcommand{\inter}{\mbox{Int}}

\newcommand{\N}{{\mathbb N}}

\newcommand{\Z}{{\mathbb Z}}
\newcommand{\R}{{\mathbb R}}
\newcommand{\Q}{{\mathbb Q}}

\DeclareMathOperator{\arsinh}{arsinh}

\DeclareMathOperator{\dist}{dist}
\DeclareMathOperator{\diam}{diam}

\DeclareMathOperator{\ir}{Int}
\DeclareMathOperator{\ao}{\alpha_1}

\DeclareMathOperator{\sg}{sgn}

\newtheorem{lemma}{Lemma} [section]
\newtheorem{thm}[lemma]{Theorem}
\newtheorem{prop}[lemma]{Proposition}
\newtheorem{cor}[lemma]{Corollary}

\newtheorem{lem}[lemma]{Lemma}
\newtheorem{rk}[lemma]{Remark}
\newtheorem{question}[lemma]{Question}
\newtheorem{conjecture}[lemma]{Conjecture}
\newtheorem{cla}[lemma]{Claim}

\title{On the intersection form of surfaces}

\author{Daniel Massart and Bjoern Muetzel}
\date{\today}

\begin{document}
\begin{abstract}
Given a closed, oriented surface $M$, the algebraic intersection of closed curves induces a symplectic form $\inter(.,.)$ on the first homology group of $M$. If $M$ is equipped with a Riemannian metric $g$, the  first homology group of $M$ inherits a norm, called the stable norm. We study the norm of the bilinear form $\inter(.,.)$, with respect to the stable norm.
\end{abstract}

\maketitle

\section{Introduction}      \label{sec:intro}  
Let $M$ be a closed (i.e. compact, without boundary) manifold of dimension two,  different from the 2-sphere, equipped with an orientation 2-form $\Omega$.  If $\alpha$ and $\beta$ are two $C^1$  closed curves on $M$ which intersect transversally, we call algebraic intersection of $\alpha $ and $\beta$  the number
$$
\sum\frac{\Omega (\dot\alpha_x, \dot\beta_x)}{| \Omega (\dot\alpha_x, \dot\beta_x) |},
$$
where \begin{itemize}
  \item $\dot\alpha_x$ denotes the tangent vector to $\alpha$ at $x$
  \item the sum is taken over all pairs of parameter values $(s,t)$ such that $\alpha (s) = \beta (t)$.
\end{itemize}
It is classical that this number only depends on the homology classes of $\alpha$ and $\beta$. We denote it by $\inter ( \left[\alpha \right],\left[\beta\right])$. The map $\inter(.,.)$ extends by linearity to a symplectic (i.e bilinear, antisymmetric, nondegenerate) form on the first homology $H_1(M,\R)$ of $M$.

The central question in this paper is \emph{when $M$ is endowed with a Riemannian metric $g$, how much can two curves of a given length intersect ?}

This amounts to evaluating the norm of the bilinear form $\inter(.,.)$ with respect to a certain norm on $H_1(M,\R)$, called the stable norm.
Informally speaking the stable norm measures the size, relative to the metric $g$, of a homology or cohomology class. Various equivalent definitions exist, see \cite{nonor, Federer74, GLP, Massart97 }. We shall use that of \cite{GLP}: for $x\in M$ and a vector $v\in T_x M$, we denote by $|v|$ its  Riemannian  norm. The {\em comass} of a differential one-form on $M$ is given by
\begin{equation}           \label{comass}           
\mbox{comass}(\omega) = \sup \big\{ \frac{|\omega(v)|}{|v|} \co x \in M, \  v \in T_x M, \ v \neq 0 \big\}.
\end{equation}

Equation~(\ref{comass}) defines a  norm on the space
$\mathcal{F}_1(M)$ of smooth 1-forms on $M$.
We get a norm on the first cohomology of $M$ by taking the infimum of the comass over all smooth closed 1-forms in a given cohomology class:
$$
\forall c \in H^1 (M,\R),\  \|c \|_s := \inf \{ \mbox{comass}(\omega) \co \omega \in \mathcal{F}_1(M), d\omega = 0, [\omega]=c \}.
$$
The norm  $ \|. \|_s$ is called the stable norm on $H^1(M,\R)$.
We denote in the same way the dual norm on $H_1(M,\R)$.

We say a homology class $h\in H_1(M,\R)$ is integer  if $h$ is the image in $ H_1(M,\R)$ of an element of $H_1(M,\Z)$.
When  $M$ is an orientable surface of genus $s$, and the homology class $h$ is integer, the stable norm has a nice expression, see \cite{nonor,Massart96, Massart97}:  $\|h\|_s$ is the minimum of all sums
$\sum |r_i| l_g (\gamma_i)$,  where
\begin{itemize}
    \item
    the index $i$ ranges over $0, \ldots, s$
    \item
    $l_g$ denotes the length with respect to $g$
    \item
    the $r_i$ are integer numbers
    \item
    the $\gamma_i$ are pairwise disjoint simple closed geodesics
    \item
    $h = \sum^{n}_{i=1} r_i \left[ \gamma_i\right]$.
\end{itemize}
The norm of the bilinear form $\inter (.,.) $ with respect to the stable norm on $H_1 (M,\R)$ is then defined as:
\begin{equation}\label{def K}
K(M,g) :=\sup \left\{ \frac{|\inter (h_1, h_2)|}{\|h_1\|_s \, \|h_2\|_s} \co  h_1,\, h_2 \in H_1(M,\R) \backslash\{0\} \right\}.
\end{equation}
Observe that in the above expression, the supremum is actually a maximum, since the function $|\inter (h_1, h_2)| / \|h_1\|_s \, \|h_2\|_s$ is zero-homogeneous, so it is actually defined on the projectivized of $H_1(M,\R)$, which is compact.

When there is no ambiguity on $M$ and $g$, we shall sometimes abbreviate the notation $K(M,g)$ to $K$.

While, from a geometrical standpoint,  the stable norm is the most natural norm on $H_1(M,\R)$, from the complex analysis viewpoint, the most natural norm is the $L^2$-norm. For any differential one-form $\omega$ and for $x \in M$ we denote by $\|\omega_x \|$ the norm, with respect to the metric $g$,  of the corresponding linear form on $T_xM$. Then we define the $L^2$-norm of $\omega$ by the formula
$$
(\| \omega \|_{2})^2  :=  \int_M \|\omega_x\|^2 d\,\vol (x),
$$
where $\vol$ denotes the volume element of the metric $g$. We define the $L^2$-norm of a cohomology class $c$ as $\inf ( \| \omega \|_{2})$,  over all 1-forms $\omega \in c$. It is a remarkable fact (see \cite{GH})  that this infimum is actually a minimum, and is achieved by the unique harmonic 1-form in the cohomology class $c$. The norm on $H_1 (M,\R)$ dual to the $L^2$-norm on $H^1(M,\R)$ will also be called $L^2$-norm, and will be denoted by the same symbol. 

The original motivation for this article was to compare the stable norm and the $L^2$-norm. This is done in Section \ref{comparison}, and our result is:
\begin{thm}        \label{compa1}              
Let \begin{itemize}
  \item $(M,g)$ be a closed, oriented surface equipped with a Riemannian metric
  \item $\emvol (M,g)$ be the total volume of $(M,g)$.
\end{itemize}
  Then for all $h \in H_1(M,\R)$, we have
\begin{equation}            \label{ineq1}    
\frac{1}{\sqrt{\mbox{\em \vol} (M,g)}} \| h \|_s \leq \|h\|_2
\leq K(M,g) \sqrt{\mbox{\em \vol} (M,g)} \| h \|_s.
\end{equation}
\end{thm}
This theorem was originally proved as Equation (4.8) of \cite{Massart96}, see also \cite{Hebda, Bangert-Katz}. In Section \ref{section compa} we give a short and simple proof.  The first inequality, which is a straightforward consequence of the Cauchy-Schwarz inequality,  has been extended to higher dimensions in \cite{Paternain}. It is also used in \cite{Haettel,Osuna}.

Now that we've been introduced to the number $K(M,g)$, we want to know more about it.
A trivial, but nice observation, is that Theorem \ref{compa1} entails
\begin{equation} \label{poincare-vol-eq}  
K(M,g) \geq  \frac{1}{\vol (M,g) }.
\end{equation}
The first question that comes to mind is
\begin{question}\label{question 1}
Is the lower bound of Equation (\ref{poincare-vol-eq}) optimal ? If not, what is the best possible lower bound ? Is it realized by some surfaces, and if so, how to characterize such surfaces ?
\end{question}
Such as it is, Question \ref{question 1} is readily answered by \cite{Hebda, Bangert-Katz}, and the answer is that $\mbox{Vol}(M,g) K(M,g) = 1$ if and only if $M$ is the two-torus and the metric $g$ is flat. The "if" part may be checked by elementary calculations and we leave it as an exercise. The "only if" holds because, by \cite{Hebda}, if the stable norm and the $L^2$-norm are proportional, then  each harmonic 1-form has constant norm. Then Proposition 6.2 of \cite{Bangert-Katz} implies that $(M,g)$ is a flat torus.

This answer to Question \ref{question 1} prompts new questions:
\begin{question} \label{question 2}
If we fix a genus $s >1$ for $M$, what is the optimal lower bound ? Is it realized by some Riemannian metrics? If so, are those metrics of constant curvature ?
\end{question}
Another obvious question is
\begin{question}\label{question 3}
Does $K(M,g)$ have an upper bound involving known geometric quantities such as the length of a homological systole (the length of a shortest, non-separating closed geodesic) ?
\end{question}
The best we can do about Questions \ref{question 2}, \ref{question 3} is summed up in Corollary \ref{bounds in arbitrary curvature} which we restate here for the commodity of the reader:
\begin{thm}
Let
\begin{itemize}
  \item $M$ be a closed,  oriented surface of genus $s \geq 1$
  \item $g$ be a Riemannian metric on $M$
  \item $D:= \diam(M,g)$ be the  diameter of $(M,g)$
  \item $l_1 := l_1 (M,g)$ be the length of a homological systole of $(M,g)$.
\end{itemize}
Then we have
$$
 \frac{1}{2l_1D} \le K(M,g) \le \frac{9}{l_1^2}.
$$
\end{thm}
In Section \ref{section constant curvature} we specialize to metrics of constant negative curvature, and we obtain the
\begin{thm}\label{encadrement courbure -1}
Let
\begin{itemize}
  \item $M$ be a closed, oriented surface of genus $s >1$
  \item $g$ be a Riemannian metric of constant curvature $-1$ on $M$
  \item $l_1$ be the length of a homological systole of $(M,g)$.
\end{itemize}
Then there exist positive numbers  $A(s)$ and $B(s)$, which depend only on the genus $s$ of $M$, such that when $l_1$ is small enough,
$$
\frac{A(s)}{l_1 |\log (l_1) |} \leq  K(M,g) \leq \frac{B(s)}{l_1 |\log (l_1) |} .
$$
\end{thm}
It would be interesting to know if there is a more precise asymptotic estimate for the behaviour of $K(M,g)$ when $g$ tends to infinity in the moduli space  $\mathcal{M}_s$ of surfaces of genus $s$ and curvature $-1$. At least we know that $K(M,g)$ does not have a maximum in $\mathcal{M}_s$, but the following question remains:
\begin{question}
Does $K(M,g)$ have a minimum when $(M,g)$ ranges over the moduli space $\mathcal{M}_s$ of surfaces of genus $s$ ? If so, which surfaces realize the minimum ?
\end{question}

There is still an obvious question that we haven't addressed:
\begin{question}
Given a surface $(M,g)$, by which homology classes is $K(M,g)$ realized, as the maximum in Equation (\ref{def K}) ? When is it realized by (the homology classes of) simple closed geodesics ?
\end{question}
In the case of flat tori, it can  be checked by elementary calculations  that for almost every  flat torus (with respect to Lebesgue measure  on the moduli space of flat tori), $K(M,g)$ is not realized by the homology classes of simple closed geodesics. In the case of surfaces of constant negative curvature, we propose the following conjecture, inspired by Theorem 10.7 of \cite{Thurston}:
\begin{conjecture}
For any $s>1$, for almost every $(M,g)$ in $\mathcal{M}_s$, $K(M,g)$ is  realized by the homology classes of simple closed geodesics.
\end{conjecture}

\section{Comparison between the stable norm and the $L^2$-norm}\label{section compa}
\subsection{Poincar\'e duality}
First let us recall some basic facts.
Let $\omega$ and $\omega'$ be two closed $1$-forms on $M$, and let $c$ and $c'$ be their respective cohomology classes. The wedge product
$\omega' \wedge \omega $ is a $2-$form on $M$, so there exists $\lambda$ in $\R$ such that
$$
\left[\omega' \wedge \omega \right] = \lambda \left[\Omega\right],
$$
 where $\Omega $ is the volume form of $M$. The number $\lambda$ only depends on the cohomology classes $c$ and $c'$, we denote it
 $c' \wedge c$ for the sake of brevity.

Recall that the Poincar\'e duality $P$ is the map from $H^1 (M,\R)$ to $H_1 (M,\R)$ induced by the wedge product of $1$-forms: for any $c, c' \in H^1 (M,\R)$, we have
\begin{equation} \label{Poincare bracket}
\langle c', P(c) \rangle = c' \wedge c,
\end{equation}
where $\langle . , . \rangle $ denotes the usual duality bracket beween $H^1 (M,\R)$ and $H_1 (M,\R)$. By \cite{GH}, p. 59, the Poincar\'e duality maps the wedge product in $H^1 (M,\R)$ to the intersection pairing in  $H_1 (M,\R)$, that is, for any $c, c' \in H^1 (M,\R)$, we have
$$
 c' \wedge c = \inter \left(P(c'), P(c) \right).
$$
We shall need the following three lemmas:
\begin{lemma} \label{Poincare L2}
The Poincar\'e duality map
\[
P: H^1(M,\R) \longrightarrow H_1(M,\R)
\]
is an isometry with respect to the $L^2$-norm.
\end{lemma}
\proof
Let $c$ be any cohomology class in $H^1(M,\R)$. Recall that the $L^2$-norm of $c$ is $\sqrt{|c \wedge ^*c|}$, where $ ^*$ is Hodge's star operator.
Then by the definition of the dual $L^2$-norm on $H_1(M,\R)$, we have
$$
\|P(c)\|_2 = \sup_{c' \neq 0} \frac{\langle c', P(c) \rangle }{\|c'\|_2}=  \sup_{c' \neq 0} \frac{c' \wedge c }{\sqrt{|c' \wedge ^*c' |}}
$$
and by the Schwarz inequality, the supremum above is achieved for $c' = ^* c$, so we have
$$
\|P(c)\|_2 = \sqrt{|c \wedge ^*c|} = \|c\|_2.
$$
\qed

\begin{lemma} \label{poincare-eq}
For any $c\in H^1(M,\R)$ and any $h\in H_1(M,\R)$ we
have
\begin{equation*}  
\langle c,h \rangle  = \inter(Pc, h).
\end{equation*}
\end{lemma}
\proof
Since the wedge product and the intersection pairing are Poincar\'e dual, we have $ \inter(Pc, h)= c \wedge P^{-1}h$, and by Equation (\ref{Poincare bracket}) we have
$$
c \wedge P^{-1}h = \langle c, P ( P^{-1}h) \rangle = \langle c,h \rangle.
$$
\qed

\begin{lemma} \label{Poincare norme K}
The norm of the inverse Poincar\'e duality map from $H_1 (M,\R)$ to $H^1 (M,\R)$, with respect to the stable norm,  is $K$.
\end{lemma}
\proof
The norm of the inverse Poincar\'e duality map from $H_1 (M,\R)$ to $H^1 (M,\R)$ is
$$
\sup_{h \neq 0} \frac{\|P^{-1}h\|_s}{\|h\|_s} = \sup_{h \neq 0} \sup_{h' \neq 0}  \frac{\langle P^{-1}h, h' \rangle }{\|h\|_s \|h'\|_s}
$$
\and by Lemma \ref{poincare-eq}
$$
\sup_{h \neq 0} \frac{\|P^{-1}h\|_s}{\|h\|_s}  =  \sup_{h \neq 0} \sup_{h' \neq 0}  \frac{ \inter(h, h') }{\|h\|_s \|h'\|_s} = K.
$$
\qed

\subsection{Proof of Theorem 1.1.} \label{comparison}

Let $\omega$ be a differential 1-form on $M$. We have
$$
(\| \omega \|_{2})^2  =  \int_M \|\omega_x\|^2 d\,\vol (x)
\leq
\int_M  (\sup_{x\in M} \|\omega_x \|^2) d\,\vol (x)  =
\vol (M,g)(\mbox{comass}(\omega))^2.
$$
Thus,
\begin{equation}            \label{L2comass}  
\| \omega \|_{2} \le \sqrt{\vol (M,g)}\ \mbox{comass}(\omega).
\end{equation}
Taking the infimum of either member of Equation (\ref{L2comass})  over all closed 1-forms in the cohomology class $c:= \left[\omega\right]$, we get
\begin{equation}  \label{thm 1, eq 1}
\| c \|_{2} \le \sqrt{\vol (M,g)} \|c\|_s
\end{equation}
and since
$$
\forall h \in H_1 (M,\R), \  \|h\|_s = \sup_{c \neq 0}\frac{\langle c, h \rangle}{\|c\|_s}
$$
we get
$$
\forall h \in H_1 (M,\R), \  \|h\|_s \leq  \sup_{c \neq 0}\sqrt{\vol (M,g)}\frac{\langle c, h \rangle}{\|c\|_2} = \sqrt{\vol (M,g)} {\|h\|_2},
$$
which is the first inequality of Theorem \ref{compa1}.

Now let us prove the second inequality. Equation (\ref{thm 1, eq 1}) and Lemma  \ref{Poincare norme K} tell us that for any homology class $h$, we have
$$
\|Ph\|_2\le \sqrt{\vol (M,g)}\ \|Ph\|_s \le \sqrt{\vol (M,g)}\ K\|h\|_s.
$$
Now Lemma \ref{Poincare L2} says that $\|Ph\|_2 = \|h\|_2$, which completes the proof.
\qed

\section{More on $K$} \label{more on K}

We get a better understanding of $K$ by noticing that in its definition we may restrict to simple closed geodesics:
\begin{lemma}\label{K simple closed geodesic}
We have
\[
K=\sup \left\{ \frac{|\inter (\alpha, \beta)|}{l_g(\alpha)l_g(\beta)} \co  \alpha, \beta \mbox{ are simple closed geodesics} \right\}.
\]
\end{lemma}
\proof
First let us point out that for any simple closed geodesics $\alpha$ and $ \beta$, we have $l_g(\alpha) \geq \| \left[\alpha\right] \|_s $ and
$l_g(\beta) \geq \| \left[\beta\right] \|_s $, so
$$
\frac{|\inter (\alpha, \beta)|}{l_g(\alpha)l_g(\beta)}  \leq  \frac{|\inter (\left[\alpha\right], \left[\beta\right])|}{\| \left[\alpha\right] \|_s \| \left[\beta\right] \|_s }\leq K,
$$
whence
\[
K \geq \sup \left\{ \frac{|\inter (\alpha, \beta)|}{l_g(\alpha)l_g(\beta)} \co  \alpha, \beta \mbox{ are simple closed geodesics} \right\}.
\]
To establish the reverse inequality, first observe that in the definition of $K$ we may restrict to integer homology classes:
\[
K=\sup \left\{ \frac{|\inter (h_1,h_2 )|}{\|h_1\|_s \|h_2\|_s} \co  h_1, h_2 \in H_1(M,\Z) \backslash\{0\} \right\}
\]
because rational homology classes are dense in $H_1(M,\R)$. Now take
\begin{itemize}
  \item two integer homology classes $h_1$ and $h_2$
  \item simple closed geodesics $\alpha_1, \ldots \alpha_k$ and $\beta_1, \ldots \beta_p$
  \item integers $a_1, \ldots a_k$ and $b_1, \ldots b_p$
\end{itemize}
such that \begin{itemize}
  \item $h_1 = a_1 \left[\alpha_1\right]+ \ldots a_k \left[\alpha_k\right]$
  \item $\|h_1\|_s = |a_1| l_g(\alpha_1 ) + \ldots |a_k| l_g (\alpha_k)$
  \item $h_2 = b_1 \left[\beta_1\right]+ \ldots b_p \left[\beta_p\right]$
  \item $\|h_2\|_s = |b_1| l_g(\beta_1 ) + \ldots |b_p| l_g (\beta_p)$.
\end{itemize}
Without loss of generality we may assume that
$$
\frac{|\inter (\alpha_1, \beta_1)|}{l_g(\alpha_1)l_g(\beta_1)}  = \max_{i,j}\frac{|\inter (\alpha_i, \beta_j)|}{l_g(\alpha_i)l_g(\beta_j)} .
$$
Then
\begin{eqnarray*}
|\inter (h_1, h_2)| = | \sum_{i,j} a_i b_j \inter (\alpha_i, \beta_j)| & \leq & \sum_{i,j} |a_i| |b_j| |\inter (\alpha_i, \beta_j)| \\
& \leq & \sum_{i,j} |a_i| |b_j| l_g(\alpha_i)l_g(\beta_j)  \frac{|\inter (\alpha_1, \beta_1)|}{l_g(\alpha_1)l_g(\beta_1)}  \\
&=& \|h_1\|_s \|h_2\|_s  \frac{|\inter (\alpha_1, \beta_1)|}{l_g(\alpha_1)l_g(\beta_1)} ,
\end{eqnarray*}
whence
$$
\frac{|\inter (h_1, h_2)|}{\|h_1\|_s \|h_2\|_s}  \leq  \frac{|\inter (\alpha_1, \beta_1)|}{l_g(\alpha_1)l_g(\beta_1)},
$$
which proves that
\[
K \leq \sup \left\{ \frac{|\inter (\alpha, \beta)|}{l_g(\alpha)l_g(\beta)} \co  \alpha, \beta \mbox{ are simple closed geodesics} \right\},
\]
and the lemma.
\qed   \\

Let $\epsilon > 0$ be a positive real number and let $\alpha$ and $\beta$ be two simple closed geodesics, such that
\[
K \leq  \frac{|\ir(\alpha,\beta)|}{l_g(\alpha) \cdot l_g(\beta)} + \epsilon .
\]

Set $N := |\ir(\alpha,\beta)|$. Replacing, if necessary, $\alpha$ and $\beta$ by shorter curves whose algebraic intersection is $N$, we may assume that $\alpha$ and $\beta$ minimize the product $l_g(\alpha) \cdot l_g(\beta)$ among all pairs of curves whose algebraic intersection is $N$. Then the following lemma tells us that all intersections of $\alpha$ and $\beta$ have the same sign, that is, the algebraic intersection of $\alpha$ and $\beta$ coincides with the number of their intersection points. Therefore in the definition of $K$ we may restrict to simple closed geodesics, where $\#\{\alpha \cap \beta\} = |\ir(\alpha,\beta)|$.

\begin{lem} Let $N$ be a positive integer, and let  $\alpha$, $\beta$  be  simple closed geodesics such that  $ |\ir(\alpha,\beta)|= N$ and 
  \begin{eqnarray*}
     \frac{|\ir(\alpha,\beta)|}{ l_g(\beta)} &=& \sup\left\{ \frac{N}{ l_g(\gamma)}  \co \gamma \text{ simple closed geodesic }, |\ir(\alpha, \gamma)|=N \right\}, \\
      \frac{|\ir(\alpha,\beta)|}{ l_g(\alpha)} &=& \sup\left\{ \frac{N}{ l_g(\gamma)}  \co \gamma \text{ simple closed geodesic }, |\ir(\beta, \gamma)|=N \right\}.
\end{eqnarray*}

Then 
\[
    N = |\ir(\alpha,\beta)|= \#\{\alpha \cap \beta\}.
\]
\label{thm:samesign}
\end{lem}

\begin{figure}[h!]
\SetLabels
\L(.38*.30) $\beta$\\
\L(.38*.51) $\alpha$\\
\L(.48*.81) $\beta_1$\\
\L(.50*.51) $\alpha_1$\\
\L(.43*.44) $p_1$\\
\L(.55*.44) $p_2$\\
\endSetLabels
\AffixLabels{%
\centerline{%
\includegraphics[height=5cm,width=5cm]{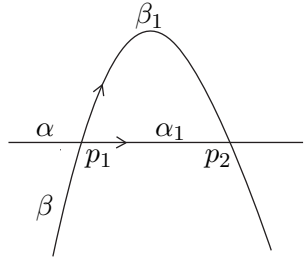}}}
\caption{Two intersections of $\alpha$ and $\beta$ with different sign.}
\label{fig:onesign}
\end{figure}

\proof  By contradiction: assume that their exists a geodesic arc $\alpha_1$ of $\alpha$ with endpoints $p_1$ and $p_2$ on $\beta$, such that the sign of the intersection at $p_1$ is different from the sign at $p_2$. Let $\beta_1$ be the geodesic arc on $\beta$ connecting $p_1$ and $p_2$ traversing $\beta$ in the positive sense and let $\beta_2$ the remaining part of $\beta$ (see Fig.~\ref{fig:onesign}).

Now assume without loss of generality that $l_g(\alpha_1) \leq l_g(\beta_1)$. We construct a new curve $\beta'$ by connecting $\beta_2$ with $\alpha_1$. By homotoping $\beta'$ away from $\alpha$ with a small deformation, the intersection points $p_1$ and $p_2$ disappear. Now the  closed geodesic $\beta''$ in the free homotopy class of $\beta'$ is strictly smaller than $\beta$ but $|\ir(\alpha,\beta)| = |\ir(\alpha,\beta'')|$. This contradicts the maximality of $|\ir(\alpha,\beta)| /  l_g(\beta)$. \qed

\begin{prop}  Let $(M,g)$ be a closed, oriented Riemannian surface. Let $l_1$ be the length of a homological systole $\ao$ and let $D$ be the diameter of $(M,g)$. Then
\begin{equation}
   \frac{1}{l_1 \cdot 2D}   \leq   K.
\label{K_low_l1}   
\end{equation}
\end{prop}
\proof
Let  $\alpha_2$ be a   shortest  closed geodesic such that $ |\ir(\ao, \alpha_2)|=1$, and let $l_2$ be its length. Then we have $K \geq (l_1 l_2)^{-1}$. We shall prove, by contradiction, that $l_2 \leq 2D$, which entails the proposition. 

Assume $l_2 > 2D$. Then there exist two points $p_1$ and $p_2$ on $\alpha_2$ whose distance is not realized by a geodesic arc on $\alpha_2$. Let $\delta_0$ be a geodesic arc, which is not an arc of $\alpha_2$, and  realizes the distance between $p_1$ and $p_2$.  The points $p_1$ and $p_2$ divide $\alpha_2$ into two arcs, $\delta_1$ and $\delta_2$. Denote by $\delta_1 \delta_0$ and $\delta_2 \delta_0$ the curves  we obtain by connecting the arcs with the same name. Both of these curves are strictly shorter than $\alpha_2$. Now we distinguish two cases.

\textit{Case 1:} $\delta_0$ and $\ao$ do not have two consecutive (along $\delta_0$) intersections with the same sign (this includes the case when $\delta_0$ and 
$\ao$ have one, or zero, intersection point). Then the algebraic intersection between $\delta_0$ and $\ao$  has absolute value zero or one. Therefore one of the curves $\delta_1 \delta_0$ and $\delta_2 \delta_0$ has algebraic intersection $\pm 1$ with $\ao$. Since it is shorter than $\alpha_2$, this contradicts the minimality of $\alpha_2$.

\begin{figure}[h!]
\SetLabels
\L(.68*.40) $\delta_0$\\
\L(.33*.62) $\ao$\\
\L(.34*.37) $p_3$\\
\L(.61*.43) $p_4$\\
\L(.45*.63) $\alpha_3$\\
\endSetLabels
\AffixLabels{%
\centerline{%
\includegraphics[height=4cm,width=8cm]{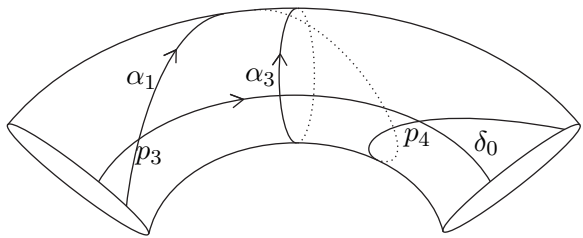}}}
\caption{The arc $\delta_0$ intersected by $\ao$ in $p_3$ and $p_4$.}
\label{fig:l1D}
\end{figure}

\textit{Case 2:} $\delta_0$ and $\ao$ have two  intersection points $p_3$ and $p_4$, consecutive along $\delta_0$,  with the same sign. Let $\alpha_3$ be the closed curve obtained by joining the arc of $\delta_0$ between $p_3$ and $p_4$, with an arc of $\ao$, of length $\leq l_1 /2$ (see Fig.~\ref{fig:l1D}). Then 
$$
l_g(\alpha_3 ) \leq l_g(\delta_0 ) + \frac{l_1}{2} \leq D +  \frac{l_1}{2} <  \frac{l_2}{2} +  \frac{l_1}{2} \leq l_2.
$$ 
On the other hand, $\alpha_3$ is homotopic to a closed curve which intersects $\ao$ exactly once. Therefore it is homotopic to a closed geodesic $c$, such that $|\ir(c,\alpha_1)|=|\ir(\alpha_3,\alpha_1)|=1$ (see \cite{dfn}, proof of Theorem 17.3.1), which is shorter than $\alpha_2$, a contradiction. 
\qed\\

\begin{figure}[h!]
\SetLabels
\L(.36*.51) $\alpha$\\
\L(.28*.60) $\beta$\\
\L(.51*.67) $I$\\
\L(.44*.90) $J$\\
\L(.44*.66) $p$\\
\L(.55*.66) $q$\\
\L(.46*.82) $\gamma'$\\
\L(.54*.84) $\gamma$\\
\endSetLabels
\AffixLabels{%
\centerline{%
\includegraphics[height=6cm,width=8cm]{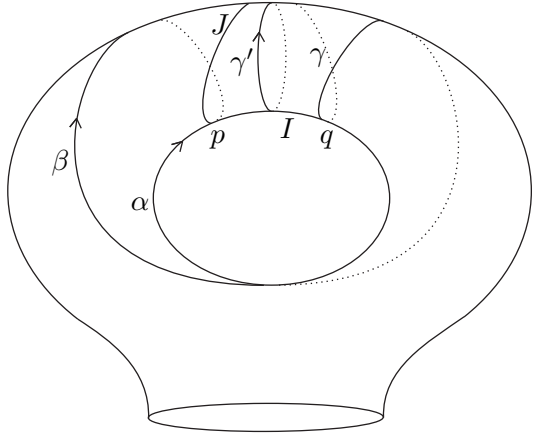}}}
\caption{A segment $I$ of $\alpha$ and a segment $J$ of $\beta$ intersecting in $p$ and $q$.}
\label{fig:intsign}
\end{figure}

In the following we will obtain an upper bound on $K$. To this end we prove the following proposition.
\begin{prop}             \label{arbitrcurvature}   
Let \begin{itemize}
  \item $(M,g)$ be a closed, oriented Riemannian surface
  \item $l_1$ be the length of a homological systole of $(M,g)$
  \item $\alpha,\beta$ be  simple closed geodesics in $(M,g)$.
\end{itemize}
Then
\begin{equation} \label{arbitrcurvature-eq}      
\frac{|\inter(\alpha,\beta)|}{l_g(\alpha)l_g(\beta)} \leq \frac{9}{l_1^{2}}.
\end{equation}

\end{prop}
\proof
Take a real number $r < l_1/2$. We cut $\alpha$ and $\beta$  into segments of length $r$ and at most one segment of smaller length.
Let $n_{\alpha}$ and $n_{\beta}$ be the respective numbers of  segments obtained. Let $I$ and $J$ be a pair of these segments in $\alpha$ and $\beta$ respectively. We shall prove, by contradiction, that the algebraic intersection of $I$ and $J$ is at most one. Assume to the contrary. Then there exist two intersection points $p$ and $q$ of $I$ and $J$, consecutive along $I$, such that the intersections of $I$ and $J$ at $p$ and $q$ have the same sign.  Let $\gamma$ be the closed curve formed by subsegments of $I$ and $J$ glued at $p$ and $q$. Then $\gamma$ is homotopic to a curve $\gamma'$ which intersects $\alpha$ exactly once (see Fig. \ref{fig:intsign}).

In particular $\gamma'$ is non-separating. On the other hand, the length of $\gamma'$ is $\leq 2r < l_1$, which contradicts the definition of $l_1$.
We have proven the inequality
\begin{equation}   \label{eq:upper1}  
|\inter(\alpha,\beta)| \leq n_{\alpha}n_{\beta}.
\end{equation}
By construction
$$
(n_{\alpha}-1)r \leq l_g(\alpha) \leq n_{\alpha}r \text{ \ \ and \ \ }
$$
$$
(n_{\beta}-1)r \leq l_g(\beta) \leq n_{\beta}r
$$
so
$$
\frac{n_{\alpha}n_{\beta}}{l_g(\alpha)l_g(\beta)}
\leq   \left(\frac{1}{r}+\frac{1}{l_g (\alpha)}\right)  \left(\frac{1}{r}+\frac{1}{l_g(\beta)}\right) \leq \left(\frac{1}{r}+\frac{1}{l_1}\right)^2.
$$
Substituting this into Equation (\ref{eq:upper1}),
and since $r$ is arbitrarily close to $l_1 /2$, we obtain the claim.
\qed

\vspace{3mm}

Summarizing Equations~(\ref{poincare-vol-eq}, \ref{K_low_l1} and \ref{arbitrcurvature-eq}),
we obtain the following bounds.

\begin{cor}     \label{bounds in arbitrary curvature}     
Let \begin{itemize}
  \item $(M,g)$ be a closed, oriented Riemannian surface
  \item $D:= \diam(M,g)$ be its diameter
  \item $l_1$ be the length of a homological systole of $(M,g)$
  \item $V:= \mbox{Vol}(M,g)$ be the total volume of $(M,g)$
  \item $K:= K(M,g)$.
\end{itemize} Then
\begin{equation*}
\label{eq:two_bounds}  
\frac{1}{V} \le K \text{ \ \ and \ \ } \frac{1}{2l_1D} \le K \le \frac{9}{l_1^2}.
\end{equation*}
\label{all_bounds}
\end{cor}

\section{Surfaces of constant negative curvature: statement of Theorem 4.2}\label{section constant curvature}
From now on we assume that the genus of $M$ is $\geq 2$ and the metric $g$ has curvature $-1$ everywhere. Recall that we denote by $\mathcal{M}_s$ the moduli space of surfaces of genus $s$, that is, the set of metrics of curvature $-1$ on $M$, modulo isometries. \\
\\
Let $\eta$ be a simple closed geodesic on $M$. Let $\omega_{\eta}$ be the supremum of all $w$, such that the geodesic arcs of length $w$ emanating  perpendicularly from  $\eta$ are pairwise disjoint. A \textit{collar} around $\eta$ or \textit{cylinder} of width $w<\omega_{\eta}$, $C_w(\eta)$, is defined by
\[
C_w(\eta):=\left\{p \in M \mid \dist(p,\eta) < w \right\}.
\]
By \cite{bu}, Theorem 4.3.2 we have:

\begin{thm}[Collar theorem]
Let $M$ be a closed, oriented surface of genus $s \geq 2$, endowed with a metric $g$ of curvature $-1$. Let $\eta$ be a simple closed geodesic in $(M,g)$.
Then
\begin{equation*}
     \omega_{\eta} \geq cl(l_g(\eta)) :=  \arsinh\left(\frac{1}{\sinh(\frac{l_g(\eta)}{2})}\right).
\end{equation*}
If $\delta$ is another simple closed geodesic that does not intersect $\eta$, then $C_{cl(l_g(\eta))}(\eta)$ and  $C_{cl(l_g(\delta))}(\delta)$ are disjoint.
\label{thm:col_lem}
\end{thm}
The main result of this section is:

\begin{thm} Let
\begin{itemize}
\item $M$ be a closed, oriented surface $M$ of genus $s \geq 2$
\item $g$ be a metric of curvature $-1$ on $M$
\item $\ao$ be a homological systole of $(M,g)$
\item $l_1 := l_g(\ao)$ be its length.
\end{itemize}
We have
\[
     \left(  (s-1)\cdot l_1 \cdot (105 s + 4\arsinh\left(\frac{4}{l_1}\right)\right) ^{-1} <  K(M,g) \leq  144 +  \frac{18(s-1)}{ l_1 \cdot cl(l_1)}.
\]
\label{thm:inequ_K}
\end{thm}

\begin{rk} Since $\arsinh\left(4/l_1\right)$ and  $cl(l_1)$ are equivalent  to $-\log (l_1)$, when $l_1$ goes to zero, Theorem \ref{thm:inequ_K} implies
 Theorem \ref{encadrement courbure -1}.
 \end{rk}
This means that $K(M,g)$ tends to infinity if  and only if $l_1$, the length of a homological systole of $(M,g)$ goes to zero.
In particular $K(M,g)$ is unbounded.
\subsection{Proof of the lower bound in Theorem 4.2.}
By \cite{bse} there exist simple closed geodesics $\beta_1=\alpha_1$ and $\beta_2, \ldots, \beta_{2s}$ such that
\begin{itemize}
  \item $\left[\beta_1\right], \ldots, \left[\beta_{2s}\right]$ is a basis of $H_1 (M, \R)$ as a vector space
  \item $\inter ( \left[\beta_{2i-1}\right], \left[\beta_{2i}\right]) = 1$ for all $i= 1, \ldots s$
  \item
  \[
 \forall i= 1, \ldots, 2s,\     l_g(\beta_i) < (s-1)\left(105 s + 4\arsinh\left(\frac{4}{l_1}\right)\right).
\]
  \end{itemize}
In particular
  \[
l_2 = l_g(\beta_2) < (s-1)\left(105 s + 4\arsinh\left(\frac{4}{l_1}\right)\right),
\]
As $\inter ( \left[\alpha_{1}\right], \left[\beta_{2}\right]) = 1$, we get:
$$
K(M,g) \geq \frac{1}{ l_1 l_2} > \frac{1}{ l_1 (s-1)\left(105 s + 4\arsinh\left(\frac{4}{l_1}\right)\right)}.
$$
\subsection{ Preliminaries to  the proof of the upper bound in Theorem 4.2.}

Let $(\alpha_i)_{i=1,..,k}$ be the set of non-separating simple closed geodesics, such that
\[
   l_g(\alpha_i) < \frac{1}{4} <  2\arsinh(1).
\]
It follows from \cite{bu}, Theorem 4.1.1 that $k \leq 3s-3$. For $i=1,..,k$ we set 
$w_i := cl(l_g(\alpha_i))-1.3$,
\[
C_i :=C_{cl(l_g(\alpha_i))}(\alpha_i) \text{ \  and  \ }  B_i :=C_{w_i}(\alpha_i) \subset C_i \text{ (see Fig. \ref{fig:collar}). }
\]
Let furthermore $\partial_1 B_i$ and $\partial_2 B_i$ be the connected boundary components  of $B_i$ and let $b_i$ be a geodesic arc realizing the distance between these two boundaries. We will gather some useful facts about the geometry of these collars in the following:
\begin{figure}[h!]
\SetLabels
\L(.77*.50) $C_i$\\
\L(.58*.65) $B_i$\\
\L(.28*.78) $C_i \backslash B_i$\\
\L(.30*.19) $\partial_1 B_i$\\
\L(.64*.19) $\partial_2 B_i$\\
\L(.49*.33) $\alpha_i$\\
\L(.42*.46) $b_i$\\
\endSetLabels
\AffixLabels{%
\centerline{%
\includegraphics[height=6cm,width=8cm,]{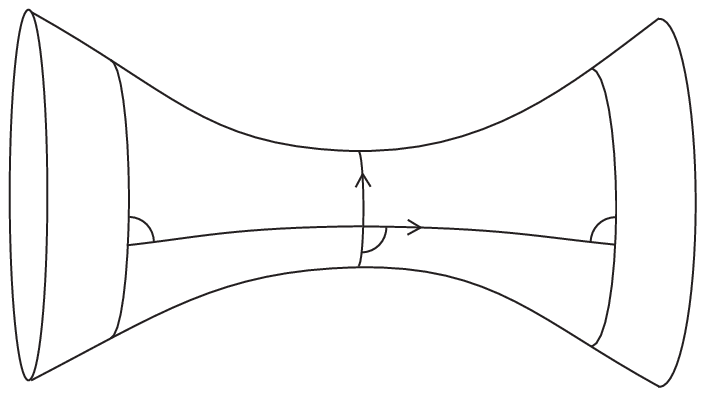}}}
\caption{An embedded collar $C_i =C_{cl(l_g(\alpha_i))}(\alpha_i)$ of a short non-separating simple closed geodesic $\alpha_i$.}
\label{fig:collar}
\end{figure}

\begin{figure}[h!]
\SetLabels
\L(.82*.02) $l_g(\alpha_i)$\\
\endSetLabels
\AffixLabels{%
\centerline{%
\includegraphics[height=8cm,width=8cm,angle=-90]{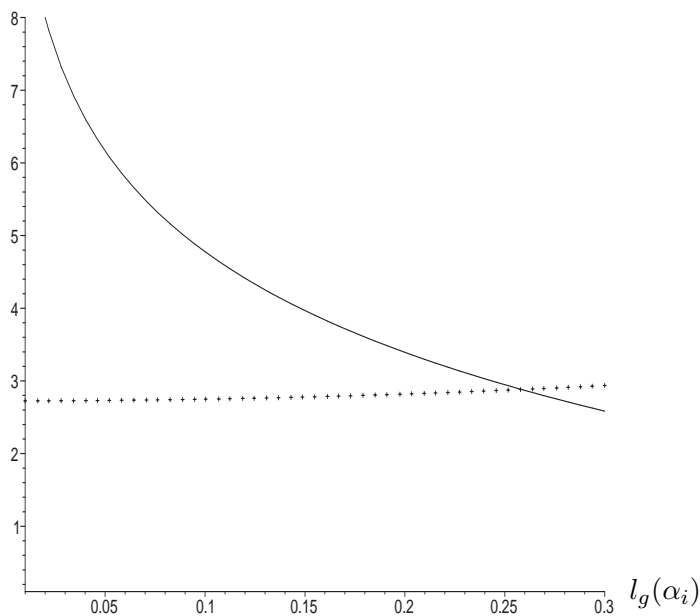}}}
\caption{Plot of $2\cdot \left(cl(l_g(\alpha_i))-1.3\right)$ \textit{(solid line)} and $5l_g(\partial_1 B_i)$ \textit{(dotted line)} in the interval $l_g(\alpha_i) \in [0.01,0.3]$.}
\label{fig:cla}
\end{figure}

\begin{itemize}
\item
By the collar theorem the $\left(C_i\right)_{i=1,..,k}$ are pairwise disjoint.
\item
For $x \leq 2\arsinh(1)$ we have: $cl(x)$ is a monotonically decreasing function and $x \leq 2 cl(x)$.
\item
$ l_g(\partial_1 B_i) = l_g(\alpha_i) \cdot \cosh(w_i)= l_g(\alpha_i) \cdot \cosh(cl(l_g(\alpha_i))-1.3).$
\item
As $l_g(\alpha_i) < \frac{1}{4}$ it follows from this formula and the definition of $cl(l_g(\alpha_i))$ that
\begin{eqnarray}
\nonumber
l_g(b_i) = 2\cdot \left(cl(l_g(\alpha_i))-1.3 \right) &>&  5 l_g(\partial_1 B_i) \text{ \ and \ } \\
l_g(\partial_1 B_i) &>& \frac{1}{2} > 2 l_g(\alpha_i). \text{ \ \ \ (see Fig. \ref{fig:cla}) \ }
\label{eq:cla}
\end{eqnarray}
\end{itemize}

Let $\epsilon > 0$ be a positive real number and let $\gamma$ and $\delta$ be two simple closed geodesics, such that
\[
K(M,g) \leq  \frac{|\ir([\gamma],[\delta])|}{l_g(\gamma) \cdot l_g(\delta)} + \epsilon .
\]
Set $N := |\ir([\gamma],[\delta])|$. Replacing, if necessary, $\gamma$ and $\delta$ by shorter curves whose algebraic intersection is $N$, we may assume that $\gamma$ and $\delta$ minimize the product $l_g(\gamma) \cdot l_g(\delta)$ among all pairs of curves whose algebraic intersection is $N$. Then Lemma \ref{thm:samesign} tells us that all intersections of $\gamma$ and $\delta$ have the same sign, that is,  the algebraic intersection  of $\gamma$ and $\delta$ coincides with the number of their intersection points.
We will obtain our result by distinguishing two cases: either  $\gamma$ or $\delta$  is one of the $(\alpha_i)_{i=1,..,k}$,  or neither  $\gamma$ nor $\delta$  is one of the $(\alpha_i)_{i=1,..,k}$.\\
\section{Proof of Theorem 4.2., Case 1: Neither $\gamma$ nor $\delta$ is one of the $(\alpha_i)_{i=1,..,k}$}
 Set
\[
     M_1 =\bigcup \limits_{i=1,..,k} B_i \text{ \ \ and \ \ }  M_2 =  M \backslash  \bigcup \limits_{i=1,..,k} B_i.
\]
For $i \in \{1,2\}$, let $N_i :=  \#\{ p \in M_i : p \in \{\gamma \cap \delta\} \}$ be the number of intersection points of  $\gamma$ and $\delta$ in $M_i$. We have:
\begin{equation}
    K(M,g) - \epsilon  \leq  \frac{N_1 + N_2}{l_g(\gamma) \cdot l_g(\delta)} = \frac{N_1}{l_g(\gamma) \cdot l_g(\delta)}+ \frac{N_2}{l_g(\gamma) \cdot l_g(\delta)} = K_1 + K_2,
\label{eq:K12}
\end{equation}
where $K_1$ and $K_2$ are defined by
$$
K_1 :=  \frac{N_1}{l_g(\gamma) \cdot l_g(\delta)} \mbox{ and } K_2 :=  \frac{N_2}{l_g(\gamma) \cdot l_g(\delta)}.
$$
Here we use the fact that the algebraic intersection  of $\gamma$ and $\delta$ coincides with the number of their intersection points. We will establish two independent upper bounds on $K_1$ and $K_2$ to prove our theorem.\\
\subsection{Upper bound on $K_1$}
To find an upper bound on $K_1$, we establish bounds on each $B_i$. Denote by
 \begin{itemize}
  \item $N^i_1 := \#\{ p \in B_i : p \in \{\gamma \cap \delta\} \}$ the number of intersection points of $\gamma$ and $\delta$ in $B_i$
  \item $K^i_1:=  \frac{ N^i_1}{l_g(\gamma) \cdot l_g(\delta)}$.
\end{itemize} . We have:
\[
   K_1 =  \frac{N_1}{l_g(\gamma) \cdot l_g(\delta)} = \sum \limits_{i=1}^k \frac{ N^i_1}{l_g(\gamma) \cdot l_g(\delta)} = \sum \limits_{i=1}^k K^i_1.
\]
Now fix an index $i$ and let us work in the cylinder $B_i$. Let
\[
\left(\gamma_j\right)_{j=1,..,n_1} =\gamma \cap B_i \text{ \ and \ } \left(\delta_l\right)_{l=1,..,n_2} = \delta \cap B_i
\]
be the disjoint union of geodesic arcs of $\gamma$ and $\delta$, respectively, which traverse $B_i$. We have:
\[
    K^i_1 = \frac{N^i_1}{l_g(\gamma) \cdot l_g(\delta)} < \frac{\sum \limits_{j,l} \# \{\gamma_j \cap \delta_l \}}{ l_g(\gamma \cap B_i)\cdot l_g(\delta \cap B_i)} =  \frac{\sum \limits_{j,l} \# \{\gamma_j \cap \delta_l \}}{\sum \limits_{j,l} l_g(\gamma_j)\cdot l_g(\delta_l)}.
\]
Now we may assume without loss of generality that
$$ \frac{ \# \{\gamma_1 \cap \delta_1 \}}{l_g(\gamma_1)\cdot l_g(\delta_1)} = \max \limits_{j,l} \frac{ \# \{\gamma_j \cap \delta_l \}}{l_g(\gamma_j)\cdot l_g(\delta_l)}. $$
It follows that
\begin{equation}
K^i_1 <  \frac{\sum \limits_{j,l} \frac{\# \{\gamma_1 \cap \delta_1 \}}{l_g(\gamma_1)\cdot l_g(\delta_1)} \cdot l_g(\gamma_j)\cdot l_g(\delta_l)}{\sum \limits_{j,l} l_g(\gamma_j)\cdot l_g(\delta_l)} = \frac{\# \{\gamma_1 \cap \delta_1 \}}{l_g(\gamma_1)\cdot l_g(\delta_1)}.
\label{eq:Ki_gd}
\end{equation}
We now determine an upper bound on the intersection number and a lower bound on the length of $\gamma_1$ and $\delta_1$. To this end we will define the winding number of an arc traversing $B_i$ and prove two lemmas concerning the intersection of two such geodesic arcs. Then we will provide a lower bound on the length of an arc.\\

The next three lemmata will also be used in the proof of the upper bound on $K_2$, where we shall need to apply them to pairs of geodesic arcs, which are not necessarily formed by an arc of $\gamma$ and an arc of $\delta$ traversing a cylinder $B_i$.

\begin{lem} Let $c$ be a geodesic arc traversing a cylinder $B_i$. Let $d \neq c$ be another geodesic arc traversing $B_i$ or let $d$ be the simple closed geodesic $\alpha_i$.
Then $c$ and $d$ intersect under the same sign at any intersection point.
\label{thm:sign_cyl}
\end{lem}

\proof
By contradiction: Assume that there exists a geodesic arc $c^1$ of $c$ with endpoints $p_1$ and $p_2$ on $d$, consecutive along $d$, such that the sign of the intersection at $p_1$ is different from the sign at $p_2$. Let $d^1$ be the geodesic arc on $d$ connecting $p_1$ and $p_2$ in $B_i$. Consider a lift of $c^1$ to the hyperbolic plane. By abuse of notation we denote the lift of $c^1$ by the same symbol. We also denote the lifts of $p_1$ and $p_2$ on $c^1$ by the same symbols. Let $d'^1$ be the lift of $d^1$ intersecting $c_1$ at $p_1$ in the hyperbolic plane. Now, due to the topology of the cylinder, $p_2$ lies also on $d'^1$. Hence in the hyperbolic plane $p_1$ and $p_2$ are connected by the two different geodesic arcs $c^1$ and $d'^1$. These arcs belong to two different geodesics passing through $p_1$ and $p_2$. But in the hyperbolic plane there can be only one geodesic through any two distinct points (see \cite{bu}, Theorem 1.1.4), a contradiction.
\qed  \\

Let $c$ be a geodesic arc traversing $B_i$. With respect to its fixed endpoints on $\partial B_i$  $c$ is in the homotopy class
\[
        [c]= [b' \cdot a  \cdot b''].
\]
Here $b'$ and $b''$ are directed geodesic arcs that meet $\alpha_i$ perpendicularly on opposite sides of $\alpha_i$, and $a$ is a directed arc on $\alpha_i$. We define the \textit{orientation} $\sigma(c)$  of $c$ with respect to $\alpha_i$ as
$$
\sigma(c) := \Bigg\{
\begin{array}{cl}
+1&  \mbox{ if the orientation of } a \mbox{ agrees with that of } \alpha_i \\
0 & \mbox{ if } a \mbox{ is a single point} \\
-1&  \mbox{ if the orientation of } a \mbox{ disagrees with  that of } \alpha_i,
\end{array}
$$
and the \textit{winding number} of $c$  as
\[
\tilde{c} := \sigma(c)\frac{l_g(a)}{l_g(\alpha_i)}=\sigma(a)\frac{l_g(a)}{l_g(\alpha_i)}
\]
so we have $\sigma(c) = \sg(\tilde{c})$, where $\sg(\cdot)$ denotes the sign function.\\
Note that due to Lemma \ref{thm:sign_cyl} an arc $c$ has only one intersection point with $\alpha_i$. Denote by the floor function the mapping $\lfloor \cdot \rfloor:\R \rightarrow \Z$, defined by $\lfloor x \rfloor:=\max\{z \in \Z : z \leq x \}$ for all $x \in \R$. The number of intersection points of two geodesic arcs traversing a cylinder $B_i$ and the sign of their intersections are related to the winding numbers of these arcs in the following way. 

\begin{lem} Let $c$ and $d$ be two geodesic arcs traversing a cylinder $B_i$. Let $\tilde{c}$ and $\tilde{d}$ be the winding numbers of the arcs $c$ and $d$, respectively. If $c$ and $d$ intersect $\alpha_i$ under the same sign, then
\[
 \lfloor |\tilde{d} - \tilde{c}| \rfloor    \leq  \#\{ c \cap d \} \leq \lfloor |\tilde{d} - \tilde{c}| \rfloor  +1.
\]
Furthermore $\sg(\tilde{d} - \tilde{c})$ determines the sign of the intersection of $c$ and $d$ at any intersection point.\\
If $c$ and $d$ intersect $\alpha_i$ under different signs, then
\[
\lfloor |\tilde{d} + \tilde{c}| \rfloor    \leq  \#\{ c \cap d \} \leq \lfloor |\tilde{d} +\tilde{c}| \rfloor  +1.
\]
Furthermore $\sg(\tilde{d} + \tilde{c})$ determines the sign of the intersection of $c$ and $d$ at any intersection point.\\
 \label{thm:up_intgd1}
\end{lem}

\proof
We first consider the case, where $c$ and $d$ intersect $\alpha_i$ from the same side.\\
\\
\textit{Case 1: $\ir(c,\alpha_i) = \ir(d,\alpha_i)$ }\\
\\
We first treat the case where $\tilde{c}=0$.\\
By Lemma \ref{thm:sign_cyl} $c$ and $d$ always intersect under the same sign at any intersection point. Hence we have
\[
\# \{c \cap d \} = |\ir(c,d)|.
\]
If $\tilde{c}=0$, then $c$ intersects $\alpha_i$ perpendicularly. Now $d$ winds at least $\lfloor |\tilde{d}| \rfloor$ times around $B_i$ and intersects $c$ at least $\lfloor |\tilde{d}| \rfloor$ times. There might be an additional intersection point, but not more than one, that is,
\begin{equation*}
    \lfloor  |\tilde{d}| \rfloor  \leq \# \{c \cap d \} = |\ir(c,d)|  \leq \lfloor  |\tilde{d}| \rfloor  +1.
\end{equation*}
Furthermore, the sign of the intersection of $c$ and $d$ at any intersection point is determined by the orientation $\sigma(d)$ of $d$, that is,
\[
      \sigma(d) = \sg(\tilde{d}) =   \left\{
 {\begin{array}{*{20}c}
   \sg(\ir(c,d))  \\
  -\sg(\ir(c,d))  \\
\end{array}} \right.  \text{ \ if \ }   
\begin{array}{*{20}c}
    \ir(c,\alpha_i) = \ir(d,\alpha_i) = +1  \\
    \ir(c,\alpha_i) = \ir(d,\alpha_i) = -1 . \\
\end{array}   
\]
Therefore our lemma is true if $\tilde{c}=0$.\\
The intersection number of two curves only depends on the homotopy class with fixed endpoints of the curves. Let $\eta$ and $\mu$ be two curves with fixed endpoints and with homotopy classes $[\eta]= [c]$ and $[\mu]=[d]$. If  $\tilde{c}=0$ we obtain the following result:
\begin{eqnarray}
\label{eq:wind01} 
    \lfloor  |\tilde{d}| \rfloor  &\leq&  |\ir(\eta, \mu )| \leq \lfloor  |\tilde{d}| \rfloor  +1 \text{ \ and \ }  \\
 \sg(\tilde{d})&=& \left\{
 {\begin{array}{*{20}c}
   \sg(\ir(\eta,\mu))  \\
  -\sg(\ir(\eta,\mu))  \\
\end{array}} \right.  \text{ \ if \ }   
\begin{array}{*{20}c}
    \ir(\eta,\alpha_i) = \ir(\mu,\alpha_i) = +1  \\
    \ir(\eta,\alpha_i) = \ir(\mu,\alpha_i) = -1.  \\
\end{array}      
\label{eq:wind02}
\end{eqnarray}
If $\tilde{c} \neq 0$, we apply a Dehn twist to the cylinder, which we define in the following via Fermi coordinates. \\
We recall that $w_i = \frac{l_g(b_i)}{2}$. The \textit{Fermi coordinates} with base point $p_1:= \alpha_i(0)$ are an injective parametrization
\[
\psi: \R \mod \{ x \mapsto x + l_g(\alpha_i) \} \times (-w_i,w_i) \rightarrow B_i, \psi: (t,s) \mapsto \psi(t,s),
\]
such that 
\begin{itemize}
   \item $\psi(0,0) = \alpha_i(0) = p_1$ and $\psi(t,0)=\alpha_i(t)$, for all $t$
   \item  $s \mapsto \psi(t,s)$ is an arc-length parametrization of an oriented geodesic arc $b_t$ that intersects $\alpha_i$ perpendicularly in $\alpha_i(t)$
   \item $\ir(b_t,\alpha_i)=+1$. 
\end{itemize}

Let $z \in \R$ be a real number. A Dehn twist of order $z$, $\mathcal{D}_z: B_i \rightarrow B_i$ is defined in Fermi coordinates by $\mathcal{D}_z(\psi(t,s)) := \psi\left(t+ z \frac{w_i +s}{2w_i},s\right)$.\\
Let  
\[
    \mathcal{D} := \mathcal{D}_{-\ir(c,\alpha_i)\tilde{c}}
\]
be the Dehn twist of order $-\ir(c,\alpha_i)\tilde{c}$. Then the winding number $\tilde{c}'$ of the geodesic arc $c'$ in the homotopy class (with fixed extremities) of $\mathcal{D}(c)$ is $0$. The winding number $\tilde{d}'$ of the geodesic arc $d'$ in the homotopy class of  $\mathcal{D}(d)$ is $\tilde{d} - \tilde{c}$.

Since $\mathcal{D}$ is isotopic to the identity, we have $ \ir(c,d) = \ir(\mathcal{D}(c),\mathcal{D}(d))$. By (\ref{eq:wind01}) we have
$$
    \lfloor  |\tilde{d'}| \rfloor  \leq \# \{c' \cap d' \} = |\ir(c',d')|  \leq \lfloor  |\tilde{d}| \rfloor  +1
$$
   that is,
\[
 \lfloor |\tilde{d} - \tilde{c}| \rfloor    \leq  \#\{ c \cap d \} =  |\ir(c,d)| \leq \lfloor |\tilde{d} - \tilde{c}| \rfloor  +1.
\]
It follows furthermore from (\ref{eq:wind02}) that $\sg(\tilde{d}')= \sg(\tilde{d} - \tilde{c})$ determines the sign of the intersection of $d$ and $c$ at any intersection point. This completes the proof in the case, where here $c$ and $d$ intersect $\alpha_i$ from the same side. \\
We now consider the case, where $c$ and $d$ intersect $\alpha_i$ from different sides.\\
\\
\textit{Case 2: $\ir(c,\alpha_i) = - \ir(d,\alpha_i)$ }\\
\\
In this case let $c^{-1}$ be the geodesic that coincides pointwise with $c$, but which traverses $B_i$ in the opposite sense. Let $\tilde{c}^{-1}$ be the winding  number of the arc $c^{-1}$. We have that
\[
 \tilde{c}^{-1}=-\tilde{c} , \text{ therefore }  \tilde{d} - \tilde{c}^{-1} = \tilde{d} + \tilde{c}.
\] 
As  $\ir(c,d) = -\ir(c^{-1},d)$ we have that
\[ 
      \# \{c \cap d \} = \# \{c^{-1} \cap d \}  \text{ \ and \ } \sg(\ir(c,d)) = -\sg(\ir(c^{-1},d)).
\] 
As  $\ir(c^{-1},\alpha_i)) = \ir(d,\alpha_i)$, we can apply the result from \textit{Case 1}. Therefore the statement for \textit{Case 2} follows from \textit{Case 1}. This completes the proof of Lemma \ref{thm:up_intgd1}.
\qed\\

We now give two lower bounds for the length of a geodesic arc traversing a cylinder $B_i$.

\begin{lem} Let $c$ be a geodesic arc traversing a cylinder $B_i$ and let $\tilde{c}$ be its winding number. We have:
\[
     l_g(c) \geq 2\cdot \left(cl(l_g(\alpha_i))-1.3\right) \text{ \ and \ }  l_g(c) \geq  |\tilde{c}| \cdot l_g(\alpha_i).
\]
\label{thm:low_length}
\end{lem}

\begin{figure}[h!]
\SetLabels
\L(.30*.62) $T$\\
\L(.40*.63) $c^*$\\
\L(.33*.51) $c^p$\\
\L(.55*.67) $B'_i$\\
\L(.51*.30) $b'_i$\\
\L(.46*.46) $q$\\
\L(.77*.49) $\alpha''_i$\\
\endSetLabels
\AffixLabels{%
\centerline{%
\includegraphics[height=6cm,width=8cm,]{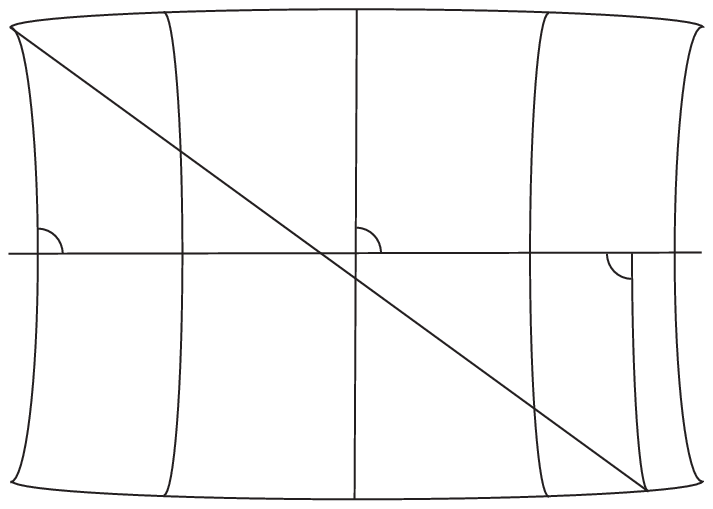}}}
\caption{Lift of the geodesic arc $c$ in the universal covering.}
\label{fig:col_lift}
\end{figure}

\proof
We lift $c$ to $c^*$ in the hyperbolic plane (see Fig.~\ref{fig:col_lift}).  Let $\alpha''_i$ be the lift of $\alpha_i$ and let $b'_i$ be a lift of $b_i$ (see Fig.~\ref{fig:collar}). Let $B'_i$ be a fundamental domain of $B_i$, whose boundary is $b'_i$.
Denote by $c^p$ the arc which we obtain from the orthogonal projection of $c^*$ onto $\alpha''_i$ in the hyperbolic plane. Let $q$ be the midpoint of $c^*$. Here, due to the symmetry of the situation, the midpoint of $c^*$ lies on $\alpha''_i$ and is also the midpoint of $c^p$. Let $T$ be a triangle with vertices $q$, an endpoint of $c^*$ and an endpoint of $c^p$ (see Fig.~\ref{fig:col_lift}).\\
It follows from the geometry of the hyperbolic right-angled triangle $T$ (see \cite{bu}, p. 454) that
\begin{equation*}
\cosh(\frac{l_g(c^*)}{2}) = \cosh(\frac{l_g(b'_i)}{2}) \cdot \cosh(\frac{l_g(c^p)}{2}).
\end{equation*}
As $l_g(b_i) = l_g(b'_i)$ and $l_g(c)=l_g(c^*)$ we obtain from the above equation that
\begin{equation*}
\cosh(\frac{l_g(c)}{2}) = \cosh(\frac{l_g(b_i)}{2}) \cdot \cosh(\frac{l_g(c^p)}{2}).
\end{equation*}
As $\cosh$ is a strictly increasing function on $\R^+$ and as $\cosh(0)=1$, it follows from the above equation that
\begin{equation*}
\cosh(\frac{l_g(c)}{2}) \geq \cosh(\frac{l_g(b_i)}{2}) \text{ \ and \ } \cosh(\frac{l_g(c)}{2}) \geq \cosh(\frac{l_g(c^p)}{2}).
\end{equation*}
From these two inequalities we obtain again by the monotonicity of the $\cosh$ function on $\R^+$ that
\begin{eqnarray*}
l_g(c) \geq l_g(b_i) = 2\cdot \left(cl(l_g(\alpha_i))-1.3\right) \text{ \ \ and \ \ }  l_g(c) \geq l_g(c^p).
\end{eqnarray*}
Here the first inequality is the first inequality of the lemma. It follows furthermore from the definition of the winding number that
\[
      |\tilde{c}| \cdot l_g(\alpha_i) = l_g(c^p) \leq l_g(c),
\]
which proves the second inequality in Lemma \ref{thm:low_length}.
\qed  \\

We will denote in the following the winding number of an arc $\gamma_j$ by $\tilde{c}_j$ and the winding number of an arc $\delta_l$ by $\tilde{d}_l$. Now let $m_1, m'_1 \in \N$ be the natural numbers such that
\[
     m_1 =  \lfloor  |\tilde{c}_1| \rfloor +1     \text{ \ \ and \ \ }  m'_1  = \lfloor |\tilde{d}_1| \rfloor  +1.
\]
We may assume without loss of generality that $m'_1 \leq m_1$. It follows from Lemma \ref{thm:up_intgd1} that
\[
   \# \{\gamma_1 \cap \delta_1 \} \leq m'_1 + m_1 \leq 2 m_1.
\]
It follows from the definition of $m_1$ and Lemma \ref{thm:low_length} that
\begin{eqnarray}
\label{eq:low_gd1}
l_g(\gamma_1) &\geq& |\tilde{c}_1|\cdot l_g(\alpha_i) \geq (m_1-1) l_g(\alpha_i) \text{ \ \ and \ \ } l_g(\gamma_1) \geq l_g(b_i) > l_g(\alpha_i)  \\
l_g(\delta_1) &\geq& l_g(b_i) = 2(cl(l_g(\alpha_i))-1.3).
\end{eqnarray}
Here the inequality $l_g(b_i) > l_g(\alpha_i)$ in (\ref{eq:low_gd1}) follows by combining the two inequalities in (\ref{eq:cla}). Using the above inequalities we obtain from inequality (\ref{eq:Ki_gd}) for $m_1 \geq 2$:
\begin{equation*}
K^i_1  \leq \frac{\#\{ \gamma_1 \cap \delta_1 \}}{l_g(\gamma_1)\cdot l_g(\delta_1)}  \leq \frac{2m_1}{(m_1-1) \cdot l_g(\alpha_i)\cdot l_g(\delta_1)} \leq \frac{2}{  l_g(\alpha_i)\left(cl(l_g(\alpha_i))-1.3\right)}.
\label{eq:Ki}
\end{equation*}
If $m_1 = 1$, we use the second inequality in (\ref{eq:low_gd1}) to derive the same upper bound.\\
\\
Combining the estimates for the $K^i_1$ we obtain that
\begin{equation*}
    K_1 \leq \sum \limits_{i=1}^k K^i_1 \leq \sum \limits_{i=1}^k \frac{2}{  l_g(\alpha_i)\left(cl(l_g(\alpha_i))-1.3\right)} <  \sum \limits_{i=1}^k  \frac{6}{l_g(\alpha_i)\cdot cl(l_g(\alpha_i))}.
\end{equation*}

\begin{figure}[h!]
\SetLabels
\L(.82*.02) $x$\\
\endSetLabels
\AffixLabels{%
\centerline{%
\includegraphics[height=8cm,width=8cm,angle=-90]{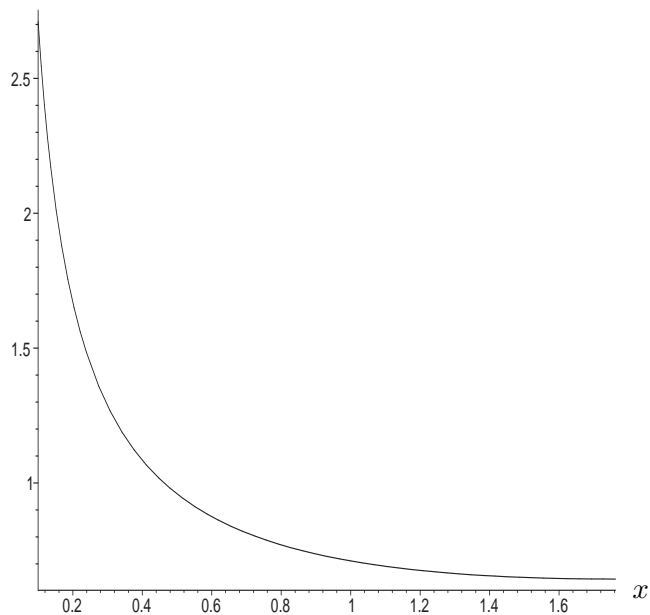}}}
\caption{Plot of the function $\frac{1}{x\cdot cl(x)}$ in the interval $[0.1,2\arsinh(1)]$.}
\label{fig:acla}
\end{figure}

To prove the last inequality we have to show that
\[
  cl(l_g(\alpha_i))-1.3 > \frac{cl(l_g(\alpha_i))}{3}  \mbox{ that is, } cl(l_g(\alpha_i)) > \frac{3}{2}\cdot 1.3 .
\]
We obtain this result by combining the two inequalities $2\cdot \left(cl(l_g(\alpha_i))-1.3 \right) >  5 l_g(\partial_1 B_i)$ and $l_g(\partial_1 B_i) > \frac{1}{2}$ in (\ref{eq:cla}).\\
As the function $\frac{1}{ x\cdot cl(x)}$ is monotonously decreasing in the interval $(0,2\arsinh(1)]$ (see Fig.~\ref{fig:acla}), we obtain from the above inequality for $K_1$ that
\begin{equation}
    K_1 \leq \sum \limits_{i=1}^k \frac{6}{  l_g(\alpha_i)\cdot cl(l_g(\alpha_i))} \leq \sum \limits_{i=1}^k \frac{6}{  l_g(\alpha_1)\cdot cl(l_g(\alpha_1))} \leq \frac{18s-18}{  l_g(\alpha_1)\cdot cl(l_g(\alpha_1))}.
\label{eq:K1}
\end{equation}
Here the last inequality follows from the fact that we have at most $3s-3$ cylinders $B_i$.

\subsection{Upper bound on $K_2$}
We recall that
\[
     K_2 = \frac{N_2}{l_g(\gamma) \cdot l_g(\delta)}, \text{ \ where \ } N_2= \#\{ p \in  M_2 = M \backslash \bigcup \limits_{i=1,..,k} B_i  : p \in \{\gamma \cap \delta\} \}.
\]

To obtain an upper bound on $K_2$ we now construct a comparison surface $M'$ from $M$ such that
\begin{equation}
l_g(\alpha'_1) > \frac{1}{4},
\label{eq:Lprime}
\end{equation}
where $\alpha'_1$ is the shortest non-separating simple closed geodesic in $M'$. Then we construct two comparison curves $\gamma' \subset M'$ and $\delta' \subset M'$, such that
\[
   l_g(\gamma') \leq  l_g(\gamma), \text{ \ \ } l_g(\delta') \leq  l_g(\delta), \text{ \ \ and \ \ } N_2 \leq |\ir([\gamma'],[\delta'])|.
\]
Now for every $\epsilon >0$, we can approximate our non-smooth surface $(M',g')$ with a smooth surface $(M_{\epsilon},g_\epsilon)$ such that the distance function of $(M_{\epsilon},g_\epsilon)$ is $\epsilon$-close to that of $(M',g')$. It follows from this remark and by applying Proposition \ref{arbitrcurvature} to $(M_{\epsilon},g_\epsilon)$ that
\begin{equation}
      K_2 = \frac{N_2}{l_g(\gamma) \cdot l_g(\delta)} \leq  \frac{|\ir([\gamma'],[\delta'])|}{l_g(\gamma')\cdot   l_g(\delta')}  \leq K(M',g') \leq \frac{9}{l_g(\alpha'_1)^2} \leq 144.
\label{eq:K2}
\end{equation}

\subsubsection{ Construction of $(M',g')$}
We construct a surface  $(M',g')$ with a singular Riemannian metric in the following way. We cut out all collars $B_i$ from $M$ and then reconnect the open ends. Here we identify the sides in the following way. For all $i \in \{1,..,k\}$, let $J_i$ be an isotopy 
\[
J_i : B_i\times[0,1] \rightarrow B_i, \text{ \ \ such that \ \ }  J_i(\cdot,0)=id  \text{ \ \  and  \ \ } J_i(B_i,1)=\partial_2 B_i.
\]
Here $J_i$ should satisfy the following condition. For all $t \in [0,1]$ and all $p_1 \in \partial_1 B_i$ 
\[
 J_i(p_1,t) \in b_{p_1} \subset B_i, \text{ \ where \ }
\]
$b_{p_1}$ is a geodesic arc in $B_i$ with endpoint $p_1$ and that intersects $\alpha_i$ perpendicularly. We define
\[
     M' := M \backslash \bigcup \limits_{i=1,..,k} B_i \mod \{ J_i(p_1,0) = J_i(p_1,1) ,\text{  for all  } p_1 \in \partial_1B_i, i \in \{1,..,k\} \}.
\]
We call $\partial B_i$ the image of $\partial_1 B_i$ in $M'$. \\
Now we have to show that the length of a non-separating simple closed curve $\eta$ in $M'$ is bigger than $\frac{1}{4}$. We distinguish two cases: either $\eta$ intersects $\bigcup \limits_{i=1,..,k}  \partial B_i$, or not.\\
Consider first the case, where $\eta$ does not intersect a $\partial B_i$ in $M'$. Now the $(\alpha_i)_{i=1,..,k}$ are the non-separating simple closed geodesics in $M$, such that $l_g(\alpha_i) < \frac{1}{4}$. In $M$ all non-separating simple closed curves of length smaller than $\frac{1}{4}$ are contained in the union $\bigcup \limits_{i=1,..,k} B_i$ of cylinders $B_i$. As $\eta \subset M$ does not intersect this set, we have that
\[
                  l_g(\eta) > \frac{1}{4}.
\]
Any simple closed curve $\eta$ in $M'$ that intersects a $\partial B_i$ either intersects a boundary of $C_i \backslash B_i$ or is contained in $C_i \backslash B_i$. In the first case $\eta$ is longer than the distance between a boundary of $C_i$ and $\partial B_i$. In the second case $\eta$ is either contractible, or is freely homotopic to $\partial B_i$, which is the shortest curve in its free homotopy class in $M'$. As the $(B_i)_{i=1,..,k}$ are chosen in a way such that for all  $i \in \{1,..,k\}$
\[
     \dist(\partial_1 C_i,\partial B_i) = 1.3 \text{ \ and \ } l_g(\partial B_i) > \frac{1}{2}  \text{ \ (see (\ref{eq:cla})), \ }
\]
we have in any case that
\[
l_g(\eta) \geq \frac{1}{2}.
\]
Summarizing these cases we obtain that the length of any non-separating simple closed curve in $M'$ is bigger than $\frac{1}{4}$. Therefore inequality (\ref{eq:Lprime}) holds. This proves our upper bound on $K_2$ in (\ref{eq:K2}).
\subsubsection{ Construction of $\gamma'$ and $\delta'$}
First we construct two comparison curves $\gamma' \subset M'$ and $\delta' \subset M'$, such that
\[
   l_g(\gamma') \leq  l_g(\gamma), \text{ \ \ } l_g(\delta') \leq  l_g(\delta) \text{ \ \ and \ \ } N_2 \leq |\ir([\gamma'],[\delta'])|.
\]
To this end we will replace all arcs of $\gamma$ and $\delta$ traversing a cylinder $C_i$ with shorter arcs in $C_i \backslash B_i$. Proceeding this way with all $\left( C_i \right)_{i=1,..,k}$, we obtain the comparison curves $\gamma'$ and $\delta'$ from $\gamma$ and $\delta$, respectively.\\


Let in the following $C_i$ be a fixed cylinder. Before we present the construction, we  will first gather some information about the way $\gamma$ and $\delta$ intersect in $B_i \subset C_i$. We recall that
\[
\left(\gamma_j\right)_{j=1,..,n_1} =\gamma \cap B_i \text{ \ and \ } \left(\delta_l\right)_{l=1,..,n_2} = \delta \cap B_i
\]
are the arcs of $\gamma$ and $\delta$ traversing $B_i$. The following lemma shows that all arcs of $\gamma$ and all arcs of $\delta$ intersect $\alpha_i$ under the same sign. More precisely:
\begin{lem}
Let \begin{itemize}
  \item $\gamma_m$ and $\gamma_j$ be two distinct arcs of $\gamma$ traversing $B_i$
  \item $\delta_k$ and $\delta_l$ be two distinct arcs of  $\delta$ traversing $B_i$.
 \end{itemize}
Then
\[
  \ir(\gamma_m,\alpha_i) = \ir(\gamma_j,\alpha_i) \text{ \ and \ }  \ir(\delta_k,\alpha_i) = \ir(\delta_l,\alpha_i).
\]

\label{thm:sameint}
\end{lem}

\begin{figure}[h!]
\SetLabels
\L(.58*.65) $B_i$\\
\L(.48*.33) $\alpha_i$\\
\L(.40*.51) $\gamma_1$\\
\L(.58*.42) $\gamma_2$\\
\L(.45*.43) $\eta$\\
\L(.34*.44) $\eta_1$\\
\L(.67*.44) $\eta_2$\\
\endSetLabels
\AffixLabels{%
\centerline{%
\includegraphics[height=6cm,width=8cm,]{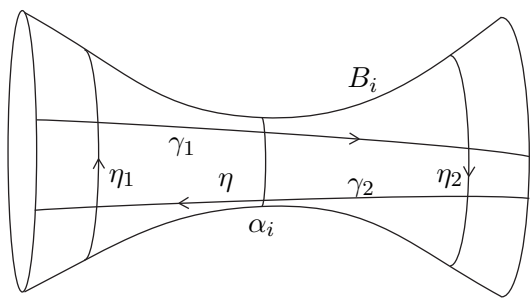}}}
\caption{A collar $B_i$ of the simple closed geodesic $\alpha_i$ and two arcs $\gamma_1$ and $\gamma_2$ of $\gamma$ traversing $\alpha_i$ under a different sign.}
\label{fig:col_g1g2}
\end{figure}

\proof
We will prove the statement by contradiction. Recall that we assume that $\gamma$ is the shortest simple closed geodesic, such that $|\ir([\gamma],[\delta])| = N$. Consider the arcs $\left(\gamma_j\right)_{j=1,..,n_1}$ of $\gamma$ traversing $B_i$. We assume that there exist two such arcs that intersect $\alpha_i$ under a different sign. Let without loss of generality  $\gamma_1$ and $\gamma_2$ be these two arcs. As $\gamma$ is a simple closed geodesic, it has no self-intersection and these two arcs can not intersect. We can therefore connect the endpoints of $\gamma_1$ and $\gamma_2$ on $\partial_1 B_i$ with an arc $\eta_1$ on $\partial_1 B_i$  and the endpoints of these two arcs on $\partial_2 B_i$ with an arc $\eta_2$ on $\partial_2 B_i$, such that together with the arcs $\gamma_1$ and $\gamma_2$ the arcs form a simple closed curve $\eta$ (see Fig. \ref{fig:col_g1g2}). As $\eta$ is contractible to a point, $[\eta] =0$. Hence replacing $\gamma_1$ and $\gamma_2$ of $\gamma$ with directed arcs corresponding to $\eta_1$ and $\eta_2$, but with inverse direction, we can create a new curve $\tilde{\gamma}$, such that $[\gamma] = [\tilde{\gamma}]$. But it follows from the first inequality in Lemma \ref{thm:low_length} and (\ref{eq:cla}) that
\[
     l_g(\gamma_j) \geq 2(cl(l_g(\alpha_i))-1.3) > l_g(\partial_1 B_i) \geq l_g(\eta_j) \text{ \ \ for \ \ } j \in \{1,2\}.
\]
Hence $l_g(\tilde{\gamma}) < l_g(\gamma)$. A contradiction to the minimality of $l_g(\gamma)$. As the same proof applies to $\delta$, we obtain our lemma.
\qed
\begin{lem}
Let $\tilde{c}_j$ be the winding number of the arc $\gamma_j$ and $\tilde{d}_j$ be the winding number of the arc $\delta_j$. Then
\[
     |\tilde{c}_j - \tilde{c}_l | < 1  \text{ \ \ and \ \ }   |\tilde{d}_j - \tilde{d}_l | < 1.
\]
\label{thm:samewind}
\end{lem}
\proof We will prove the lemma for two arcs $\gamma_j$ and $\gamma_l$ of $\gamma$. If $|\tilde{c}_j - \tilde{c}_l | \geq 1$, then $\lfloor |\tilde{c}_j - \tilde{c}_l | \rfloor \geq 1$. It follows from Lemma \ref{thm:up_intgd1} that $\# \{\gamma_j \cap \gamma_l\} \geq 1$ and $\gamma_j$ intersects $\gamma_l$. A contradiction to the fact that $\gamma$ is a simple closed geodesic and therefore has no self-intersection.
\qed \\

Let $\gamma_{min}$ and $\delta_{min}$ be two arcs of $\gamma$ and $\delta$, respectively, with minimal absolute value of the winding number. Let $\tilde{c}_{min}$ and $\tilde{d}_{min}$ be the winding numbers of these arcs. We have
\begin{equation*}
     |\tilde{c}_{min}|    = \min \limits_j |\tilde{c}_j| \text{ \ and \  } |\tilde{d}_{min}|    = \min \limits_l |\tilde{d}_l|.
\end{equation*}
We set furthermore
\begin{equation*}
\lfloor |\tilde{c}_{min}| \rfloor    = m_\gamma     \text{ \ \ and \ \ } \lfloor |\tilde{d}_{min}| \rfloor  = m_\delta.
\end{equation*}
Let $c$ be a geodesic arc traversing $B_i$. With respect to its fixed endpoints on $\partial B_i$  $c$ is in the homotopy class $[c]= [b' \cdot a  \cdot b'']$. Here $b'$ and $b''$ are directed geodesic arcs that meet $\alpha_i$ perpendicularly on opposite sides of $\alpha_i$, and $a$ is a directed arc on $\alpha_i$. We recall that $\sigma(c)$,  the orientation of $c$, is defined by
\[
\sigma(c) := \Bigg\{
\begin{array}{cl}
+1&  \mbox{ if the orientation of } a \mbox{ agrees with that of } \alpha_i \\
0 & \mbox{ if } a \mbox{ is a single point} \\
-1&  \mbox{ if the orientation of } a \mbox{ disagrees with  that of } \alpha_i.
\end{array}
\]
If $\tilde{c}$ is the winding number of the arc $c$, then $\sigma(c)= \sg(\tilde{c})$. We have:
\begin{lem}
Let $\tilde{c}_j$ be the winding number of the arc $\gamma_j$ and $\tilde{d}_j$ be the winding number of the arc $\delta_j$. \\
If $|\tilde{c}_{min} | \geq 1 $ then $\sigma(\gamma_j) = \sigma(\gamma_{min})$ for all $j \in \{1,..,n_1\}$.\\
If $|\tilde{d}_{min} | \geq 1 $ then $\sigma(\delta_l) = \sigma(\delta_{min})$ for all $l \in \{1,..,n_2\}$.\\
Furthermore
\[
 |\tilde{c}_j - \sigma(\gamma_j)m_\gamma | < 2  \text{ \ \ and \ \ }   |\tilde{d}_j - \sigma(\delta_l)m_\delta | < 2.
\]
\label{thm:less2}
\end{lem}
\proof
We will prove the lemma for an arc $\gamma_j$ of $\gamma$. The statement about the orientation $\sigma(\gamma_j)$ of $\gamma_j$ follows from Lemma \ref{thm:samewind}.
It follows from the triangle inequality that
\begin{eqnarray*}
     |\tilde{c}_j - \sigma(\gamma_j)m_\gamma | =    |\tilde{c}_j - \tilde{c}_{min} + \tilde{c}_{min}-\sigma(\gamma_j)m_\gamma | &\leq&  \\
     |\tilde{c}_j - \tilde{c}_{min}| + | \tilde{c}_{min}-\sigma(\gamma_j)m_\gamma | &<& 1 +  | \tilde{c}_{min}-\sigma(\gamma_j)m_\gamma |.
\end{eqnarray*}
Here the last inequality follows from Lemma \ref{thm:samewind}. Now if $|\tilde{c}_{min} | < 1$ then $m_\gamma = 0$ and the inequality is true. If $|\tilde{c}_{min} | \geq 1 $ then $\sigma(\gamma_j) = \sigma(\gamma_{min})$ and the inequality follows from the definition of  $ |\tilde{c}_{min}| $ and $m_\gamma$.
\qed \\

We now define the comparison curves $\gamma'$ and $\delta'$ of $M'$. Let $v$ be a geodesic arc of $\gamma$ or $\delta$ that traverses $B_i$ with endpoints $p_1 \in \partial_1 B_i$ and  $p_2 \in \partial_2 B_i$ and let $\tilde{v}$ be its winding number. We first replace $v$ with the geodesic arc $v'$ with the same endpoints $p_1$ and $p_2$ on $\partial B_i$, such that its winding number $\tilde{v}'$ has the following value: \\
If $|\tilde{c}_{min}| \leq |\tilde{d}_{min}|$ then  $m_\gamma \leq m_\delta$ and we set
\begin{equation}
    \tilde{v}' = \left\{ {\begin{array}{*{20}c}
   {\tilde{v}-\sigma(v) \max \{m_\gamma -1,0 \} }  \\
   {\tilde{v}-\sigma(v) (\max \{m_\gamma -1,0\} + \max \{(m_\delta-m_\gamma -2),0\})  }  \\
\end{array}} \right. \text{ \ if \ } \begin{array}{*{20}c}
   {v \subset \gamma}  \\
   {v \subset \delta}  \\
\end{array}.
\label{eq:changewind1}
\end{equation}
If $|\tilde{d}_{min}| \leq |\tilde{c}_{min}|$ then $m_\delta \leq m_\gamma$ and we set
\begin{equation*}
    \tilde{v}' = \left\{ {\begin{array}{*{20}c}
   {\tilde{v}-\sigma(v) \max \{m_\delta -1,0\} }  \\
   {\tilde{v}-\sigma(v) (\max \{m_\delta -1,0\} + \max \{(m_\gamma-m_\delta -2),0\})
    }  \\
\end{array}} \right. \text{ \ if \ } \begin{array}{*{20}c}
   {v \subset \delta}  \\
   {v \subset \gamma}  \\
\end{array}.
\end{equation*}

Denote by $\gamma^*$ and $\delta^*$ the curves which we obtain this way from $\gamma$ and $\delta$. We call
\[
\left(\gamma^*_j\right)_{j=1,..,n_1} =\gamma^* \cap B_i \text{ \ and \ } \left(\delta^*_l\right)_{l=1,..,n_2} = \delta^* \cap B_i
\]
the arcs of $\gamma^*$ and $\delta^*$ traversing $B_i$ and denote by  $\tilde{c}^*_j$ and $\tilde{d}^*_l$ the winding number of the arc $\gamma^*_j$ and $\delta^*_l$, respectively.\\
Now let $v''$ be the arc on $\partial B_i \subset M'$, such that
\[
v'' = J_i(v',1).
\]
Replacing all arcs $v$ of $\gamma$ and $\delta$ in all $(B_i)_{i=1,..,k}$ with corresponding arcs $v''$ in $M'$ we obtain $\gamma'$ and $\delta'$.

\begin{cla}\label{claim}
 $l_g(\gamma') \leq  l_g(\gamma)$, $l_g(\delta') \leq  l_g(\delta)$ and $N_2 \leq |\ir([\gamma'],[\delta'])|.$
\label{thm:claim}
\end{cla}
Before launching into the proof of Claim \ref{claim}, let us explain the idea a little bit.

On the one hand, we need the comparison curves $\gamma'$ and $\delta'$ to be shorter than the original curves $\gamma$ and $\delta$, respectively. So we ensure that the winding numbers of the arcs of $\gamma^*$ and the arcs of $\delta^*$ in $B_i$ are no greater than the winding numbers of the corresponding arcs of $\gamma$ and $\delta$. This is the reason for the $ \max \{m_\gamma -1,0\}$ and $ \max \{m_\delta -1,0\}$ in the definition.

On the other hand, to prove the statement about the intersection number: $N_2 \leq |\ir([\gamma'],[\delta'])|$, what we need to do is to make sure that  the intersections between the comparison curves $\gamma^*$ and $\delta^*$  have the same sign as those of the original curves. Recall that by Lemma \ref{thm:samesign} $\gamma$ and $\delta$ intersect always under the same sign and have
\[
N_2=  \#\{ p \in  M_2 = M \backslash \bigcup \limits_{i=1,..,k} B_i  : p \in \{\gamma \cap \delta\} \}
\]
intersection points outside the union of the cylinders $(B_i)_{i=1,..,k}$. As $\gamma$ and $\gamma^*$ and $\delta$ and $\delta^*$ coincide in $M_2$, $\gamma^*$ and $\delta^*$ have at least $N_2$ intersection points and their sign of intersection at any intersection point in $M_2$ is the same. We will show that due to the $\max \{(m_\gamma-m_\delta -2),0\}$ in our definition there are no two consecutive intersections of $\gamma^*$ and $\delta^*$ with different sign. It follows that
$ \# \{ \gamma^* \cap \delta^*\} = |\ir(\gamma^*,\delta^*)|$, whence $N_2 \leq |\ir(\gamma^*,\delta^*)|$.

As furthermore $\gamma'$ and $\delta'$ are the image of $\gamma^*$ and $\delta^*$, respectively, under a continuous deformation of the surface $M$, it follows that
\[
\ir(\gamma^*,\delta^*) = \ir(\gamma',\delta').
\]
In total we obtain:
\[
   N_2 \leq \# \{ \gamma^* \cap \delta^*\} = |\ir(\gamma^*,\delta^*)| = |\ir(\gamma',\delta')|.
\]
Now let us prove Claim \ref{claim}. To simplify our proof we may assume without loss of generality that
\[
|\tilde{c}_{min}| \leq |\tilde{d}_{min}|.
\]
We will first show:\\
\\
\textit{ $l_g(\gamma') \leq  l_g(\gamma)$ and $l_g(\delta') \leq  l_g(\delta)$. }\\
\\
To this end we first show that $|\tilde{c}^*_j | < 3$ and $|\tilde{d}^*_j | < 5$.\\
Consider an arc $\gamma^*_j$ of $\gamma^*$. It follows from Equation (\ref{eq:changewind1}) that
\[
   \tilde{c}^*_{j} =  \tilde{c}_j - \sigma(\gamma_j)  \max \{m_\gamma -1,0 \}.
\]
From which follows by the triangle inequality that
\begin{eqnarray*}
\nonumber
|\tilde{c}^*_j | = |\tilde{c}_j - \sigma(\gamma_j) (m_\gamma  - m_\gamma + \max \{m_\gamma -1,0 \} )|  &\leq& \\
           |\tilde{c}_j - \sigma(\gamma_j) m_\gamma | +|\sigma(\gamma_j)|\cdot | m_\gamma - \max \{m_\gamma -1,0 \} | &<& 2 +1 = 3.
\end{eqnarray*}
Here the last inequality follows from Lemma \ref{thm:less2} and the fact that $| x - \max \{x -1,0 \} | \leq 1$ for all $x \in \R_+$. \\

Now let $\delta^*_j$ be an arc of $\delta^*$. It follows from Equation (\ref{eq:changewind1}) that
\[
     \tilde{d}^*_l = \tilde{d}_l- \sigma(\delta_l)(\max \{m_\gamma -1,0 \} +  \max \{(m_\delta-m_\gamma -2),0\}).
\]
We distinguish two cases, $m_\delta -m_\gamma \leq 2$ and $m_\delta -m_\gamma > 2$.\\

\begin{enumerate}[i)]
\item $m_\delta -m_\gamma \leq 2 \Rightarrow \max \{(m_\delta - m_\gamma-2),0\}=0 $\\

We obtain from the triangle inequality and Lemma \ref{thm:less2} that
\begin{eqnarray*}
  |\tilde{d}^*_l| = |\tilde{d}_l- \sigma(\delta_l)\max \{m_\gamma -1,0 \}| =  \\
  |\tilde{d}_l- \sigma(\delta_l) (m_\delta -m_\delta + \max \{m_\gamma -1,0 \})| &\leq& \\
   |\tilde{d}_l- \sigma(\delta_l) m_\delta| + |\sigma(\delta_l)|\cdot|m_\delta - \max \{m_\gamma -1,0 \}| &<&\\
    2 + |m_\delta - \max \{m_\gamma -1,0 \}|.
\end{eqnarray*}
Applying again the triangle inequality we have that
\begin{eqnarray*}
 |\tilde{d}^*_l| < 2 + |m_\delta -m_\gamma +m_\gamma -  \max \{m_\gamma -1,0 \}|  &<& \\
     2 + |m_\delta -m_\gamma| + |m_\gamma -  \max \{m_\gamma -1,0 \}| &<& 2 + 2 +1 = 5.
\end{eqnarray*}
Here the last inequality follows from the hypothesis and the fact that the function  $| x - \max \{x -1,0 \} | \leq 1$ for all $x \in \R_+$. Hence  $|\tilde{d}^*_l| <5$.\\

\item $m_\delta -m_\gamma > 2 \Rightarrow \max \{(m_\delta-m_\gamma -2),0\}= m_\delta -m_\gamma -2 $\\

It follows from Equation (\ref{eq:changewind1}) that in this case
\[
     \tilde{d}^*_l = \tilde{d}_l- \sigma(\delta_l)(\max \{m_\gamma -1,0 \}+ m_\delta - m_\gamma -2) .
\]
In this case we apply the triangle inequality twice to $| \tilde{d}^*_l|$  and obtain
\begin{eqnarray*}
 |\tilde{d}^*_l| = |\tilde{d}_l- \sigma(\delta_l)(\max \{m_\gamma -1,0 \} + m_\delta  - m_\gamma  - 2) | &\leq& \\
 |\tilde{d}_l- \sigma(\delta_l) m_\delta| +|\sigma(\delta_l)|\cdot|m_\gamma -\max \{m_\gamma -1,0 \} | + |2| &<& 5.
\end{eqnarray*}
Here the last inequality follows from Lemma \ref{thm:less2} and the fact that the function  $| x - \max \{x -1,0 \} | \leq 1$ for all $x \in \R_+$.
\end{enumerate}
Hence in any case $|\tilde{c}^*_j | < 3$ and $|\tilde{d}^*_j | < 5$.\\
Let $v$ be an arc of $\gamma$ or $\delta$ traversing $B_i$. Let $v'$ be the replacement arc of $v$ according to Equation (\ref{eq:changewind1}). The deformation $J_i$ collapses the cylinder $B_i$ onto $\partial B_i$. Here the arc $v'$ of $\gamma^*$ or $\delta^*$ is deformed into an arc $v''$ of $\gamma'$ or $\delta'$ in $M'$. If $\tilde{v}'$ is the winding number of an arc $v'$,  then $v''$ winds $|\tilde{v}'|$ times around $\partial B_i$. Hence, as $|\tilde{v}'| < 5$, we have that
\[
    l_g(v'') <  5 l_g(\partial_1 B_i).
\]
As $v$ traverses $B_i$ we conclude with Lemma \ref{thm:low_length} that
\[
      l_g(v) > 2(cl(l_g(\alpha_i))-1.3) >  5 l_g(\partial_1 B_i) >  l_g(v'').
\]
Here the second inequality follows from the first part of inequality (\ref{eq:cla}). Hence our replacement arc $v''$ of $\gamma'$ or $\delta'$ is always shorter than the arc $v$ from $\gamma$ or $\delta$. This proves the first part of Claim \ref{thm:claim}. Now we prove the second part, that is,
$$
N_2 \leq |\ir([\gamma'],[\delta'])|.
$$
This amounts to showing  that  there are no two consecutive intersections of $\gamma^*$ and $\delta^*$ with different sign. For this it suffices  to show that the sign of the intersection between $\gamma^*$ and $\delta^*$ at any point $p$ inside a cylinder $B_i$ is always equal to $\sg(\ir(\gamma,\delta))$. We will therefore show that all comparison arcs $\gamma^*_j$ and $\delta^*_l$ intersect under the same sign as $\gamma_j$ and $\delta_l$ in a cylinder $B_i$.\\
\\
Therefore we have to treat two cases. Either both $\gamma$ and $\delta$ intersect $\alpha_i$ from the same side or $\gamma$ and $\delta$ intersect $\alpha_i$ from different sides. Here we will deduce the result of the second case from the result of the first case.\\ 
\\
\textit{Case I: $\ir(\gamma,\alpha_i) = \ir(\delta,\alpha_i)$}\\
\\
By Lemma \ref{thm:up_intgd1} it is sufficient to show that the sign of the difference of the winding numbers does not change, that is,
\begin{equation}
     \sg( \tilde{d}^*_l - \tilde{c}^*_j ) =  \sg( \tilde{d}_l - \tilde{c}_j ).
\label{cond_sign}
\end{equation}
Recall that we assume that $|\tilde{c}_{min}| \leq |\tilde{d}_{min}|$. In this case we have
\begin{equation}
     \max \{m_\gamma -1,0\} \leq m_\gamma \leq m_\delta \text{ \ \ and \ \ }   \max \{m_\delta -1,0\} \leq m_\delta.
     \label{eq:max_gd}
\end{equation}

We will prove our statement depending on whether $\tilde{c}_j$  and $\tilde{d}_l$  have different sign or are equal to zero or whether they have the same sign or are equal to zero.\\
\begin{enumerate}[a)]
\item \textit{($\tilde{c}_j \leq 0$ and $\tilde{d}_l \geq 0$) or ($\tilde{c}_j \geq 0$ and $\tilde{d}_l \leq 0$)  }\\

We assume without loss of generality that $\tilde{c}_j \leq 0 $ and $\tilde{d}_l \geq 0$. We have that
\[
 \tilde{d}_l - \tilde{c}_j  \geq 0.
\]
As $\tilde{c}_j \leq 0 $ and $\tilde{d}_l \geq 0$ it follows that $ \tilde{d}_l - \tilde{c}_j  = 0$ if and only if $\tilde{d}_l = \tilde{c}_j = 0$. This case implies that $|\tilde{c}_{min}| = |\tilde{d}_{min}|=0.$  It follows from Equation (\ref{eq:changewind1}) that
$\tilde{d}^*_l = \tilde{d}_l$  and $\tilde{c}^*_j =  \tilde{c}_j$. Hence
\[
 \tilde{d}_l - \tilde{c}_j  =  \tilde{d}^*_l - \tilde{c}^*_j = 0.
\]
and (\ref{cond_sign})  holds. \\
Conversely if $\tilde{d}_l = \tilde{c}_j = 0$ does not hold then either  $\tilde{d}_l$ or $\tilde{c}_j$ is different from zero. We may assume that $\tilde{d}_l > 0$ and we have that
\begin{equation}
 \tilde{d}_l - \tilde{c}_j  > 0 \mbox{ that is, }  \sg(\tilde{d}_l - \tilde{c}_j)  =+1 .
\label{eq:sign_a}
\end{equation}
It follows from Equation (\ref{eq:changewind1}) and from inequality (\ref{eq:max_gd}) that
\[
    \tilde{c}^*_{j} = \tilde{c}_{j} + \max \{m_\gamma -1,0 \} \leq \tilde{c}_{j} + m_\gamma \leq 0.
\]
Here the last inequality follows from the fact that $\tilde{c}_j \leq 0 $ and  $|\tilde{c}_j| \geq m_\gamma$. Furthermore
\[
    \tilde{d}^*_l =\tilde{d}_l- (\max \{m_\gamma -1,0 \} +  \max \{(m_\delta-m_\gamma -2),0\} ).
\]
We now show that $\tilde{d}^*_l >0$. As $\tilde{c}^*_{j} \leq 0$ it follows that
\[
    \tilde{d}^*_l - \tilde{c}^*_l > 0 \mbox{ that is, }  \sg(\tilde{d}^*_l - \tilde{c}^*_l) =+1
\]
and condition (\ref{cond_sign}) follows with (\ref{eq:sign_a}).\\
To show that $ \tilde{d}^*_l > 0$ we distinguish the two subcases  $m_\delta -m_\gamma \leq 2$ and  $m_\delta -m_\gamma > 2$.\\

\begin{enumerate}[i)]
\item $m_\delta -m_\gamma \leq 2 \Rightarrow  \max \{(m_\delta-m_\gamma -2),0\}=0$ \\

Then as $m_\delta \geq m_\gamma$
\[
   \tilde{d}^*_l =\tilde{d}_l- \max \{m_\gamma -1,0 \} \geq \tilde{d}_l- \max \{m_\delta -1,0 \} > 0.
\]
Here the last inequality follows by distinguishing the cases $m_\delta > 1$ and $m_\delta \leq 1$. In the latter case we use the fact that $\tilde{d}_l >0$.\\

\item $m_\delta -m_\gamma > 2 \Rightarrow  \max \{(m_\delta-m_\gamma -2),0\}= m_\delta-m_\gamma -2$ \\

Then it follows with inequality (\ref{eq:max_gd}):
\begin{eqnarray*}
   \tilde{d}^*_l =\tilde{d}_l - \max \{m_\gamma -1,0 \} - m_\delta + m_\gamma  + 2 &=& \\
                   (\tilde{d}_l - m_\delta)  + (m_\gamma - \max \{m_\gamma -1,0 \})  + 2 &\geq& 2 > 0.
\end{eqnarray*}
\end{enumerate}

\item \textit{($\tilde{c}_j \leq 0$ and $\tilde{d}_l \leq 0$) or ($\tilde{c}_j \geq 0$ and $\tilde{d}_l \geq 0$)}\\

We assume without loss of generality that $\tilde{c}_j \geq 0 $ and $\tilde{d}_l \geq 0$. We will further subdivide this case into the subcases  $m_\delta - m_\gamma \leq 2$ and $m_\delta - m_\gamma >2$.\\

\begin{enumerate}[i)]
\item   $m_\delta - m_\gamma \leq 2$\\
\\
In this case it follows from Equation (\ref{eq:changewind1}) that
\begin{eqnarray*}
    \tilde{c}^*_{j} = \tilde{c}_{j} - \max \{m_\gamma -1,0 \} \text{ \ and \ }
       \tilde{d}^*_{l} = \tilde{d}_{l} - \max \{m_\gamma -1,0 \} .
\end{eqnarray*}
Hence
\[
    \tilde{d}^*_{l} - \tilde{c}^*_{j} =    \tilde{d}_{l} - \tilde{c}_{j} \Rightarrow \sg(\tilde{d}^*_l - \tilde{c}^*_j) =  \sg(\tilde{d}_l - \tilde{c}_j)
\]
and the condition (\ref{cond_sign}) is fulfilled.\\

\item   $m_\delta - m_\gamma > 2$\\

We have due to the definition of $m_\delta$ that
\begin{eqnarray*}
     \tilde{d}_l - \tilde{c}_j \geq m_\delta - \tilde{c}_j =
     m_\delta - \tilde{c}_j + |\tilde{c}_{min}| - |\tilde{c}_{min}| + m_\gamma - m_\gamma &=&\\
   m_\delta - (\tilde{c}_j - |\tilde{c}_{min}|) - (|\tilde{c}_{min}| - m_\gamma) - m_\gamma.
\end{eqnarray*}
We note that due to Lemma \ref{thm:samewind} $0 \leq \tilde{c}_j - |\tilde{c}_{min}|  < 1$. Furthermore due to the definition of $|\tilde{c}_{min}|$ and $m_\gamma$ we have that $0 \leq |\tilde{c}_{min}| - m_\gamma < 1$.
Hence
\[
     \tilde{d}_l - \tilde{c}_j  > m_\delta - 2 - m_\gamma > 0 \Rightarrow  \sg(\tilde{d}_l - \tilde{c}_j)=+1.
\]
Here the last inequality follows from our hypothesis. It follows from Equation (\ref{eq:changewind1}) that
\begin{eqnarray*}
\tilde{c}^*_{j} = \tilde{c}_{j} -  \max \{m_\gamma -1,0 \}   \text{ \ \ and \ \ } \\   \tilde{d}^*_{l} = \tilde{d}_{l} -  \max \{m_\gamma -1,0 \} - m_\delta + m_\gamma +2 .
\end{eqnarray*}
Hence
\begin{eqnarray*}
   \tilde{d}^*_{l} - \tilde{c}^*_{j}  = \tilde{d}_{l}  - m_\delta + m_\gamma +2 -  \tilde{c}_{j} &=& \\
       2 + (\tilde{d}_{l}  - m_\delta) - ( \tilde{c}_{j} - m_\gamma) &>& 2 + 0 - 2 > 0
\end{eqnarray*}
Here by the definition of $m_\delta$, $\tilde{d}_{l}  - m_\delta \geq 0$ and by Lemma \ref{thm:less2} $\tilde{c}_{j} - m_\gamma < 2$. Hence $\sg(\tilde{d}^*_{l} - \tilde{c}^*_{j}) =+1$. Therefore
\[
     \sg(\tilde{d}^*_l - \tilde{c}^*_j) =  \sg(\tilde{d}_l - \tilde{c}_j)=+1
\]
and the condition (\ref{cond_sign}) is fulfilled.
\end{enumerate}
\end{enumerate}
This proves that in any case condition (\ref{cond_sign}) is fulfilled. This settles the claim in the case, where $\gamma$ and $\delta$ intersect $\alpha_i$ from the same side.\\ 
\\
\textit{Case II: $\ir(\gamma,\alpha_i) = -\ir(\delta,\alpha_i)$}\\
\\
If $\gamma$ and $\delta$ intersect $\alpha_i$ from different sides, it is sufficient to show that 
\begin{equation}
     \sg( \tilde{d}^*_l + \tilde{c}^*_j ) =  \sg( \tilde{d}_l + \tilde{c}_j ).
\label{cond_sign2}  
\end{equation}
This follows from Lemma \ref{thm:up_intgd1}.\\
In this case let $\gamma^{-1}$ be the oriented geodesic that coincides pointwise with $\gamma$, but which has opposite orientation. Let $(\gamma^{-1}_j)_{j=1,..,n_1}$ be the arcs of $\gamma^{-1}$ traversing $B_i$. Let $\tilde{c}^{-1}_j$ be the winding  number of the arc $\gamma^{-1}_j$. We have that $\tilde{c}^{-1}_j=-\tilde{c}_j$.\\
Applying Equation (\ref{eq:changewind1}) to the arcs of $\delta$ and $\gamma^{-1}$, we obtain that
\[
 \tilde{c}^{-1}_j=-\tilde{c}_j \text{ \ and \ }  \tilde{c}^{-1*}_j = -\tilde{c}^*_j.
\] 
Now $\gamma^{-1}$ and $\delta$ intersect $\alpha_i$ under the same sign. Hence condition (\ref{cond_sign}) from \textit{Case I} is fulfilled, that is,
\[ 
      \sg(\tilde{d}_l - \tilde{c}^{-1}_j) =   \sg(\tilde{d}^*_l - \tilde{c}^{-1*}_j).  
\]
Combining the two previous equations we obtain: 
\[
   \sg(\tilde{d}_l - (-\tilde{c}_j)) =  \sg(\tilde{d}_l - \tilde{c}^{-1}_j) =   \sg(\tilde{d}^*_l - \tilde{c}^{-1*}_j) = \sg(\tilde{d}^*_l - (-\tilde{c}^{*}_j)) 
\]
and therefore condition (\ref{cond_sign2}) is fulfilled. This settles the claim in the case where $\gamma$ and $\delta$ intersect $\alpha_i$ from different sides. Hence the remaining second part of Claim \ref{thm:claim} is true. \qed  \\

This finishes the construction of the curves $\gamma'$ and $\delta'$. From this we get the  upper bound on $K_2$  in Equation (\ref{eq:K2}). 
\subsection{End of the proof of Case 1 of Theorem 4.2.}
Summarizing \textit{Case 1}, we obtain from the inequalities (\ref{eq:K1}) and (\ref{eq:K2}) in (\ref{eq:K12}) that
\[
    K(M,g)- \epsilon \leq K_1 + K_2 \leq \frac{18s-18}{l_g(\alpha_1)cl(l_g(\alpha_1))} + 144.
\]
As $\epsilon$ is arbitrarily small, we obtain the upper bound stated in Theorem \ref{thm:inequ_K} from this inequality. 
\section{Proof of Theorem 4.2., Case 2: Either $\gamma$ or $\delta$ is one of the $(\alpha_i)_{i=1,..,k}$} 

We may suppose that $\gamma=\alpha_i$. In this case we have to verify the upper bound in Theorem \ref{thm:inequ_K} for
\[
 K(M,g) - \epsilon = \frac{|\ir([\alpha_i],[\delta])|}{l_g(\alpha_i) \cdot l_g(\delta)}= \frac{N}{l_g(\alpha_i) \cdot l_g(\delta)}.
\]
Now $\delta$ intersects $\alpha_i$ $N$ times. To this end it has to traverse $N$ times the collar $C_i$ of $\alpha_i$. Analogous to Lemma \ref{thm:low_length} we obtain from the length of the $N$ arcs $\left(\delta_j\right)_{j=1,..,N}$  of $\delta$ traversing $C_i$:
\[
l_g(\delta) > \sum \limits_{j=1}^N l_g(\delta_j)  \geq N\cdot2 cl(l_g(\alpha_i)).
\]
Now from the monotonicity of the function $\frac{1}{x\cdot cl(x)}$ in the interval $(0,2\arsinh(1)]$ (see Fig. \ref{fig:acla}) it follows that
\[
     K(M,g) - \epsilon \leq \frac{N}{2N l_g(\alpha_i)\cdot cl(l_g(\alpha_i))} = \frac{1}{2 l_g(\alpha_1)\cdot cl(l_g(\alpha_1))}.
\]
Again we obtain our upper bound in Theorem \ref{thm:inequ_K} as $\epsilon$ is arbitrarily small.\qed

\section*{Acknowledgement} While working on this article the second author has been supported by the Alexander von Humboldt foundation.\\

\vspace{2cm}

\noindent Daniel Massart and Bjoern Muetzel\\	
\noindent Department of Mathematics, Universit\'e Montpellier 2 \\
\noindent place Eug\`ene Bataillon, 34095 Montpellier cedex 5, France \\
\noindent E-Mail: \textit{massart@math.univ-montp2.fr} and \textit{bjorn.mutzel@gmail.com}

\end{document}